\documentstyle[11pt]{article}

\newtheorem{eg}{Example}
\newtheorem{lemma}{Lemma}[section]
\newtheorem{theo}{Theorem}
\newtheorem{prop}[lemma]{Proposition}
\newtheorem{cor}[theo]{Corollary}

\newtheorem{remark}{Remark}

\font\bbb=msbm10 scaled\magstep1

\newcommand{\RR}{\mbox{\bbb R}}

\newcommand{\ZZ}{\mbox{\bbb Z}}

\newcommand{\TPSSD}{S^{\hspace{.2mm}d-1}\! \times \hspace{-3.8mm}_{-} \, S^1}

\newcommand{\La}{\Lambda}

\newcommand{\si}{\sigma}

\begin{document}

\title{Three dimensional pseudomanifolds on eight vertices}
\author{Basudeb Datta,$^1$ and Nandini Nilakantan$^2$}

\date{}


\maketitle

\vspace{-7mm}

\noindent {\small $^1$Department of Mathematics, Indian Institute
of Science, Bangalore 560\,012,  India \newline \mbox{}
\hspace{1mm} dattab@math.iisc.ernet.in}


\noindent {\small $^2$Department of Mathematics \& Statistics,
Indian Institute of Technology, Kanpur 208\,016, India \newline
\mbox{} \hspace{1mm} nandini@iitk.ac.in}

\begin{center}
\date{To appear in `International Journal of Mathematics and Mathematical Sciences'}
\end{center}

\noindent {\bf Abstract.} {\small A normal pseudomanifold is a
pseudomanifold in which the links of simplices are also
pseudomanifolds. So, a normal 2-pseudomanifold triangulates a
connected closed 2-manifold. But, normal $d$-pseudomanifolds form
a broader class than triangulations of connected closed
$d$-manifolds for $d \geq 3$. Here, we classify all the 8-vertex
neighbourly normal 3-pseudomanifolds. This gives a classification
of all the 8-vertex normal 3-pseudomanifolds. There are 73 such
3-pseudomanifolds, 38 of which triangulate the 3-sphere and other
35 are not combinatorial 3-manifolds. These 35 triangulate six
distinct topological spaces. As a preliminary result, we show
that any 8-vertex 3-pseudomanifold is equivalent by proper
bistellar moves to an 8-vertex neighbourly 3-pseudomanifold. This
result is the best possible since there exists a $9$-vertex
non-neighbourly $3$-pseudomanifold ($B^3_9$ in Example \ref{N39}
below) which does not allow any proper bistellar moves.}

\bigskip

{\small

\noindent Mathematics Subject Classification (2000): 57Q15,
57Q05, 57N05, 55M25.

\smallskip

\noindent Keywords: Bistellar moves; Normal pseudomanifold;
Triangulated 3-sphere; Singular vertex.}



\section{Introduction}

Recall that a {\em simplicial complex} is a collection of
non-empty finite sets (sets of {\em vertices}) such that every
non-empty subset of an element is also an element. For $i \geq
0$, the elements of size $i + 1$ are called the {\em
$i$-simplices} (or {\em $i$-faces}) of the complex.

A simplicial complex is usually thought of as a prescription for
construction of a topological space by pasting geometric
simplices. The space thus obtained from a simplicial complex $K$
is called the {\em geometric carrier} of $K$ and is denoted by
$|K|$. We also say that $K$ {\em triangulates} $|K|$. A {\em
combinatorial $2$-manifold} (respectively, {\em combinatorial
$2$-sphere}) is a simplicial complex which triangulates a closed
surface (respectively, the 2-sphere $S^{\hspace{.2mm}2}$).

For a simplicial complex $K$, the maximum of $k$ such that $K$
has a $k$-simplex is called the {\em dimension} of $K$. A
$d$-dimensional simplicial complex $K$ is called {\em pure} if
each simplex of $K$ is contained in a $d$-simplex of $K$. A
$d$-simplex in a pure $d$-dimensional simplicial complex is
called a {\em facet}. A $d$-dimensional pure simplicial complex
$K$ is called a {\em weak pseudomanifold} if each $(d -
1)$-simplex of $K$ is contained in exactly two facets of $K$.

With a pure simplicial complex $K$ of dimension $d \geq 1$, we
associate a graph $\Lambda(K)$ as follows. The vertices of
$\Lambda(K)$ are the facets of $K$ and two vertices of
$\Lambda(K)$ are adjacent if the corresponding facets intersect
in a $(d-1)$-simplex of $K$. If $\Lambda(K)$ is connected then
$K$ is called {\em strongly connected}. A strongly connected weak
pseudomanifold is called a {\em pseudomanifold}. Thus, for a
$d$-pseudomanifold $K$, $\Lambda(K)$ is a connected $(d +
1)$-regular graph. This implies that $K$ has no proper subcomplex
which is also a $d$-pseudomanifold. (Or else, the facets of such
a subcomplex would provide a disconnection of $\Lambda(X)$.)

For any set $V$ with $\#(V) = d + 2$ ($d \geq 0$), let $K$ be the
simplicial complex whose simplexes are all the non-empty proper
subsets of $V$. Then $K$ is a $d$-pseudomanifold and triangulates
the $d$-sphere $S^{\hspace{.2mm}d}$. This $d$-pseudomanifold $K$
is called the {\em standard $d$-sphere} and is denoted by
$S^{d}_{d+2}(V)$ (or  $S^{d}_{d+ 2}$). By convention,
$S^{\hspace{.2mm}0}_2$ is the only $0$-pseudomanifold.

If $\si$ is a face of a simplicial complex  $K$ then the {\em
link} of $\si$ in $K$, denoted by ${\rm lk}_K(\sigma)$ (or ${\rm
lk}(\sigma)$), is by definition the simplicial complex whose
faces are the faces $\tau$ of $K$ such that $\tau$ is disjoint
from $\si$ and $\si\cup\tau$ is a face of $K$. Clearly, the link
of an $i$-face in a weak $d$-pseudomanifold is a weak $(d - i -
1)$-pseudomanifold. For $d\geq 1$, a connected weak
$d$-pseudomanifold is said to be a {\em normal
$d$-pseudomanifold} if the links of all the simplices of
dimension $\leq d - 2$ are connected. Thus, any connected
triangulated $d$-manifold (triangulation of a closed $d$-manifold)
is a normal $d$-pseudomanifold. Clearly, the normal
2-pseudomanifolds are just the connected combinatorial
$2$-manifolds. But, normal $d$-pseudomanifolds form a broader
class than connected triangulated $d$-manifolds for $d \geq 3$.

Observe that if $X$ is a normal pseudomanifold then $X$ is a
pseudomanifold. (If $\La(X)$ is not connected then, since $X$ is
connected, $\La(X)$ has two components $G_1$ and $G_2$ and two
intersecting facets $\sigma_1$, $\sigma_2$ such that $\sigma_i
\in G_i$, $i = 1, 2$. Choose $\sigma_1$, $\sigma_2$ among all
such pairs such that $\dim(\sigma_1 \cap \sigma_2)$ is maximum.
Then $\dim(\sigma_1 \cap \sigma_2) \leq d-2$ and ${\rm
lk}_X(\sigma_1 \cap \sigma_2)$ is not connected, a contradiction.)
Notice that all the links of positive dimensions (i.e., the links
of simplices of dimension $\leq d - 2$) in a normal
$d$-pseudomanifold are normal pseudomanifolds. Thus, if $K$ is a
normal 3-pseudomanifold then the link of a vertex in $K$ is a
combinatorial 2-manifold. A vertex $v$ of a normal
3-pseudomanifold $K$ is called {\em singular} if the link of $v$
in $K$ is not a 2-sphere. The set of singular vertices is denoted
by ${\rm SV}(K)$. Clearly, the space $|K| \setminus {\rm SV}(K)$
is a pl 3-manifold. If ${\rm SV}(K) = \emptyset$ (i.e., the link
of each vertex is a 2-sphere) then $K$ is called a {\em
combinatorial $3$-manifold}. A {\em combinatorial $3$-sphere} is a
combinatorial 3-manifold which triangulates the topological
3-sphere $S^{\hspace{.2mm}3}$.

Let $M$ be a weak $d$-pseudomanifold. If $\alpha$ is a $(d -
i)$-face of $M$, $0 < i \leq d $, such that ${\rm lk}_M(\alpha) =
S^{i - 1}_{i + 1}(\beta)$ and $\beta$ is not a face of $M$ (such
a face $\alpha$ is said to be a {\em removable} face of $M$) then
consider the weak $d$-pseudomanifold (denoted by
$\kappa_{\alpha}(M)$) whose facet-set is $\{\sigma : \sigma
\mbox{ a facet of } M, \alpha \not\subseteq \sigma\} \cup \{
\beta \cup \alpha \setminus \{v\} : v \in \alpha\}$. The
operation $\kappa_{\alpha} : M \mapsto \kappa_{\alpha}(M)$ is
called a {\em bistellar $i$-move}. For $0< i < d$, a bistellar
$i$-move is called a {\em proper bistellar move}. If
$\kappa_{\alpha}$ is a proper bistellar $i$-move and ${\rm
lk}_M(\alpha) = S^{i-1}_{i+1}(\beta)$ then $\beta$ is a removable
$i$-face of $\kappa_{\alpha}(M)$ (with ${\rm
lk}_{\kappa_{\alpha}(M)}(\beta) = S^{d-i - 1}_{d - i +
1}(\alpha)$) and $\kappa_{\beta} : \kappa_{\alpha}(M) \mapsto M$
is an bistellar $(d - i)$-move. For a vertex $u$, if ${\rm
lk}_M(u) = S^{d-1}_{d+1}(\beta)$ then the bistellar $d$-move
$\kappa_{\{u\}} : M \mapsto \kappa_{\{u\}}(M) = N$ deletes the
vertex $u$ (we also say that $N$ is obtained from $M$ by {\em
collapsing} the vertex $u$). The operation $\kappa_{\beta} : N
\mapsto M$ is called a bistellar {\em $0$-move} (we also say that
$M$ is obtained from $N$ by {\em starring} the vertex $u$ in the
facet $\beta$ of $N$). The 10-vertex combinatorial 3-manifold
$A^3_{10}$ in Example \ref{M310} is not neighbourly and does not
allow any bistellar 1-move. In \cite{bd6}, Bagchi and Datta have
shown that if the number of vertices in a non-neighbourly
combinatorial 3-manifold is at most 9 then the 3-manifold admits
a bistellar 1-move. Existence of the 9-vertex 3-pseudomanifold
$B^3_9$ in Example \ref{N39} shows that Bagchi and Datta's result
is not true for 9-vertex 3-pseudomanifolds. Here we prove.

\begin{theo}$\!\!\!${\bf .} \label{t1}
If $M$ is an $8$-vertex $3$-pseudomanifold, then there exist a
sequence of bistellar $1$-moves $\kappa_{A_1}, \dots,
\kappa_{A_m}$, for some $m\geq 0$, such that $\kappa_{A_m}(\cdots
(\kappa_{A_1}(M)))$ is a neighbourly $3$-pseudomanifold.
\end{theo}

In \cite{a1}, Altshuler has shown that every combinatorial
3-manifold with at most 8 vertices is a combinatorial 3-sphere.
In \cite{gs}, Gr\"{u}nbaum and Sreedharan have shown that there
are exactly 37 polytopal 3-spheres on 8 vertices (namely,
$S^3_{8, 1} \dots, S^3_{8, 37}$ in Examples \ref{m38n} and
\ref{m38nn}). They have also constructed the non-polytopal sphere
$S^3_{8, 38}$. In \cite{b1}, Barnette proved that there is only
one more non-polytopal 8-vertex 3-sphere (namely, $S^3_{8, 39}$).
In \cite{e}, Emch constructed a 8-vertex normal 3-pseudomanifold
(namely, $N_1$ in Example \ref{pm38}) as a block design. This is
not a combinatorial 3-manifold  and its automorphism group is
${\rm PGL}(2, 7)$ (cf. \cite{k3}). In \cite{a2}, Altshuler has
constructed another 8-vertex normal 3-pseudomanifold (namely,
$N_{5}$ in Example \ref{pm38}). In \cite{l1}, Lutz has shown that
there exist exactly three 8-vertex normal 3-pseudomanifolds which
are not combinatorial 3-manifolds (namely, $N_{1}$, $N_{5}$ and
$N_{6}$ in Example \ref{pm38}) with vertex-transitive
automorphism groups. Here we prove\,:

\begin{theo}$\!\!\!${\bf .} \label{t2} Let $S^3_{8,35}, \dots,
S^3_{8,38}$, $N_{1}, \dots, N_{15}$ be as in Examples $\ref{m38n}$
and $\ref{pm38}$.
\begin{enumerate} \item[$(a)$] Then $S^3_{8,i} \not\cong S^3_{8,
j}$, $N_{k} \not \cong N_{l}$ and $S^3_{8,m} \not\cong N_{n}$ for
$35 \leq i< j \leq 38$, $1 \leq k < l \leq 15$, $35 \leq m \leq
38$ and $1\leq n \leq 15$.
\item[$(b)$] If $M$ is an $8$-vertex neighbourly normal
$3$-pseudomanifold then $M$ is isomorphic to one of $S^3_{8,35},
\dots, S^3_{8,38}, N_{1}, \dots, N_{15}$.
\end{enumerate}
\end{theo}

\begin{cor}$\!\!\!${\bf .} \label{t3}
There are exactly $39$ combinatorial $3$-manifolds on $8$
vertices, all of which are combinatorial $3$-spheres.
\end{cor}

\begin{cor}$\!\!\!${\bf .} \label{t4} There are exactly $35$
normal $3$-pseudomanifolds on $8$ vertices which are not
combinatorial $3$-manifolds. These are $N_{1}, \dots, N_{35}$
defined in Examples $\ref{pm38}$ and $\ref{pm38nn}$.
\end{cor}

The topological properties of these normal $3$-pseudomanifolds
are given in Section 3.

\section{Preliminaries}

All the simplicial complexes considered in this paper are finite
(i.e., with finite vertex-set). The vertex-set of a simplicial
complex $K$ is denoted by $V(K)$. We identify the 0-faces of a
complex with the vertices. The $1$-faces of a complex $K$ are also
called the {\em edges} of $K$.

If $K$, $L$ are two simplicial complexes, then an {\em
isomorphism} from $K$ to $L$ is a bijection $\pi : V(K) \to V(L)$
such that for $\si \subseteq V(K)$, $\si$ is a face of $K$ if and
only if $\pi (\si)$ is a face of $L$. Two complexes $K$, $L$ are
called {\em isomorphic} when such an isomorphism exists. We
identify two complexes if they are isomorphic. An isomorphism
from a complex $K$ to itself is called an {\em automorphism} of
$K$. All the automorphisms of $K$ form a group under composition,
which is denoted by ${\rm Aut}(K)$.

For a face $\sigma$ in a simplicial complex $K$, the number of
vertices in ${\rm lk}_K(\sigma)$ is called the {\em degree} of
$\sigma$ in $K$ and is denoted by $\deg_K(\sigma)$ (or by
$\deg(\sigma)$). If every pair of vertices of a simplicial
complex $K$ form an edge then $K$ is called {\em neighbourly}.
For a simplicial complex $K$, if $U \subseteq V(K)$ then $K[U]$
denotes the induced complex of $K$ on the vertex-set $U$.

If the number of $i$-faces of a $d$-dimensional simplicial
complex $K$ is $f_i(K)$ ($0\leq i\leq d$), then the number
$\chi(K) := \sum_{i= 0}^{d}(- 1)^i f_i(K)$ is called the {\em
Euler characteristic} of $K$.

A {\em graph} is a simplicial complex of dimension $\leq 1$. A
finite 1-pseudomanifold is called a {\em cycle}. An $n$-cycle is
a cycle on $n$ vertices and is denoted by $C_n$ (or by $C_n(a_1,
\dots, a_n)$ if the edges are $a_1a_2, \dots, a_{n-1}a_n,
a_na_1$).

For a simplicial complex $K$, the graph consisting of the edges
and vertices of $K$ is called the {\em edge-graph} of $K$ and is
denoted by ${\rm EG}(K)$. The complement of ${\rm EG}(K)$ is
called the {\em non-edge graph} of $K$ and is denoted by ${\rm
NEG}(K)$. For a weak 3-pseudomanifold $M$ and an integer $n \geq
3$, we define the graph $G_n(M)$ as follows. The vertices of
$G_n(M)$ are the vertices of $M$. Two vertices $u$ and $v$ form
an edge in $G_n(M)$ if $uv$ is an edge of degree $n$ in $M$.
Clearly, if $M$ and $N$ are isomorphic then $G_n(M)$ and $G_n(N)$
are isomorphic for each $n$.

If $M$ is a weak 3-pseudomanifold and $\kappa_{\alpha} : M
\mapsto \kappa_{\alpha}(M) = N$ is a bistellar 1-move then, from
the definition, $(f_0(N)$, $f_1(N), f_2(N), f_3(N)) = (f_0(M),
f_1(M) + 1, f_2(M) + 2, f_3(M) + 1)$ and $\deg_{N}(v) \geq
\deg_M(v)$ for any vertex $v$. If $\kappa_{\alpha} : M \mapsto
\kappa_{\alpha}(M) = L$ is a bistellar $3$-move then $(f_0(L),
f_1(L), f_2(L), f_3(L)) = (f_0(M)-1, f_1(M) - 4, f_2(M) -6,
f_3(M) -3)$.

\setlength{\unitlength}{3mm}

\begin{picture}(45,11.5)(0,-2)


\thicklines

\put(2,5){\line(4,3){4}} \put(2,5){\line(3,-1){3}}
\put(2,5){\line(1,-1){4}} \put(8,5){\line(-2,3){2}}
\put(8,5){\line(-3,-1){3}} \put(8,5){\line(-1,-2){2}}
\put(5,4){\line(1,4){1}} \put(5,4){\line(1,-3){1}}

\put(14,5){\line(4,3){4}} \put(14,5){\line(3,-1){3}}
\put(14,5){\line(1,-1){4}} \put(20,5){\line(-2,3){2}}
\put(20,5){\line(-3,-1){3}} \put(20,5){\line(-1,-2){2}}
\put(17,4){\line(1,4){1}} \put(17,4){\line(1,-3){1}}
\put(18,1){\line(0,1){3.1}} \put(18,4.6){\line(0,1){3.4}}

\put(9,6){\vector(1,0){4}} \put(13,4){\vector(-1,0){4}}
\put(32,6){\vector(1,0){4}} \put(36,4){\vector(-1,0){4}}

\put(24,4){\line(1,1){4}} \put(24,4){\line(2,-1){4}}
\put(28,2){\line(0,1){6}} \put(31,4){\line(-3,-2){3}}
\put(31,4){\line(-3,4){3}}

\put(37,4){\line(1,1){4}} \put(37,4){\line(2,-1){4}}
\put(41,2){\line(0,1){6}} \put(44,4){\line(-3,-2){3}}
\put(44,4){\line(-3,4){3}}

\thinlines

\put(2,5){\line(1,0){2.9}} \put(8,5){\line(-1,0){2.4}}

\put(14,5){\line(1,0){2.9}} \put(20,5){\line(-1,0){1.7}}

\put(17.5,4.7){\mbox{-}}

\put(24,4){\line(1,0){3.6}} \put(28.4,4){\line(1,0){2.6}}

\put(37,4){\line(1,0){3}} \put(41.4,4){\line(1,0){2.6}}

\put(40,5){\line(1,-3){1}} \put(40,5){\line(-3,-1){3}}
\put(40,5){\line(1,3){1}} \put(44,4){\line(-4,1){2.6}}

\put(40.2,4.8){\mbox{$.$}} \put(40.65,4.752){\mbox{$.$}}
\put(40.5,3.948){\mbox{$.$}}

\put(9,7){\mbox{{\small $1$-move}}} \put(9.3,2.5){\mbox{{\small
$2$-move}}}

\put(32,7){\mbox{{\small $0$-move}}} \put(32.3,2.5){\mbox{{\small
$3$-move}}}

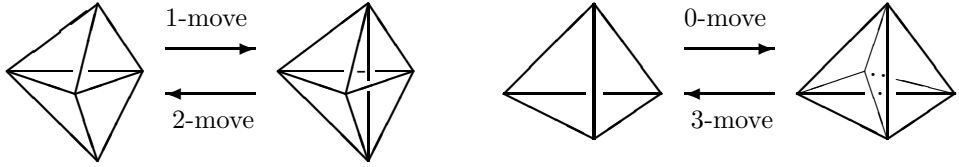
\put(14,-1){\mbox{Bistellar moves in dimension $3$}}

\end{picture}

Consider the binary relation `$\leq$' on the set of weak
3-pseudomanifolds as: $M \leq N$ if there exists a finite sequence
of bistellar $1$-moves $\kappa_{\alpha_1}, \dots,
\kappa_{\alpha_m}$, for some $m\geq 0$, such that $N =
\kappa_{\alpha_m}(\cdots \kappa_{\alpha_1}(M))$. Clearly, this
$\leq$ is a partial order relation.

Two weak $d$-pseudomanifolds $M$ and $N$ are {\em bistellar
equivalent} (denoted by $M\sim N$) if there exists a finite
sequence of bistellar operations leading from $M$ to $N$. If
there exists a finite sequence of proper bistellar operations
leading from $M$ to $N$ then we say $M$ and $N$ are {\em properly
bistellar equivalent} and we denote this by $M \approx N$.
Clearly, `$\sim$' and `$\approx$' are equivalence relations on
the set of pseudomanifolds. It is easy to see that $M \sim N$
implies that $|M|$ and $|N|$ are pl homeomorphic.

For two simplicial complexes $X$ and $Y$ with disjoint vertex
sets, the simplicial complex $X \ast Y := X \cup Y \cup
\{\sigma\cup \tau ~ : ~ \sigma \in X, ~ \tau \in Y\}$ is called
the {\em join} of $X$ and $Y$.

Let $K$ be an $n$-vertex (weak) $d$-pseudomanifold. If $u$ is a
vertex of $K$ and $v$ is not a vertex of $K$ then consider the
simplicial complex $\Sigma_{uv}K$ on the vertex set $V(K) \cup
\{v\}$ whose set of facets is $\{\sigma \cup \{u\} : \sigma$ is a
facet of $K$ and $u \not\in \sigma\} \cup \{\tau \cup \{v\} :
\tau$ is a facet of $K\}$. Then $\Sigma_{uv}K$ is a (weak)
$(d+1)$-pseudomanifold and $|\Sigma_{uv}K|$ is the topological
suspension $S|K|$ of $|K|$ (cf. \cite{bd2}). Easy to see that the
links of $u$ and $v$ in $\Sigma_{uv}K$ are isomorphic to $K$.
This $\Sigma_{uv}K$ is called the {\em one-point suspension} of
$K$.

For two $d$-pseudomanifolds $X$ and $Y$, a simplicial map $f
\colon X \to Y$ is called a {\em $k$-fold branched covering}
(with discrete branch locus) if $|f||_{|X|\setminus f^{- 1}(U)}
\colon |X| \setminus f^{- 1}(U) \to |Y| \setminus U$ is a
$k$-fold covering for some $U \subseteq V(Y)$. (We say that $X$
is a {\em branched cover} of $Y$ and $Y$ is a {\em branched
quotient} of $X$.) The smallest such $U$ (so that
$|f||_{|X|\setminus f^{- 1}(U)} \colon |X| \setminus f^{- 1}(U)
\to |Y|\setminus U$ is a covering) is called the {\em branch
locus}. If $N$ is a $k$-fold branched quotient of $M$ and
$\widetilde{N}$ is obtained from $N$ by collapsing a vertex
(respectively starring a vertex in a facet) then $\widetilde{N}$
is the branched quotient of $\widetilde{M}$, where
$\widetilde{M}$ can be obtained from $M$ by collapsing $k$
vertices (respectively starring $k$ vertices in $k$ facets). For
proper bistellar moves we have\,:

\begin{lemma}$\!\!\!${\bf .} \label{bcover}
Let $M$ and $N$ be two $d$-pseudomanifolds and $f \colon M \to N$
be a $k$-fold branched covering. For $1 \leq l < d-1$, if $\alpha$
is a removable $l$-face, then $f^{- 1}(\alpha)$ consists of $k$
removable $l$-faces $\alpha_1, \dots, \alpha_k$ $($say$)$ and
$\kappa_{\alpha_k}(\cdots (\kappa_{\alpha_1}(M)))$ is a $k$-fold
branched cover of $\kappa_{\alpha}(N)$.
\end{lemma}

\noindent {\bf Proof.} Let ${\rm lk}_N(\alpha) = S^{d-l-1}_{d-l+
1}(\beta)$. Since the dimension of $\alpha$ is $> 0$,
$f^{-1}(\alpha)$ consists of $k$ $l$-faces, $\alpha_1, \dots,
\alpha_k$ (say) of $M$. Let ${\rm lk}_M(\alpha_i) = S^{d-l-1}_{d
- l+ 1}(\beta_i)$ and $M_i := M[\alpha_i \cup \beta_i]$ for $1\leq
i \leq k$. Since $f$ is simplicial, $\beta_i$ is not a face of $M$
and hence $\alpha_i$ is removable for each $i$. Since $0< l < d -
1$, it follows that $M_i$ is neighbourly. For $i \neq j$, if $x
\neq y \in V(M_i) \cap V(M_j)$, then $xy$ is an edge in $M_i\cap
M_j$ and hence the number of edges in $f^{- 1}(f(x)f(y))$ is less
than $k$, a contradiction. So, $\#(V(M_i) \cap V(M_j)) \leq 1$ for
$i\neq j$. This implies that $\beta_i$ is not a face in
$\kappa_{\alpha_j}(M)$ and hence $\alpha_i$ is removable in
$\kappa_{\alpha_j}(M)$ for $i \neq j$. The result now follows.
\hfill $\Box$

\bigskip

Remark \ref{remark} shows that Lemma \ref{bcover} is not true for
$l = d-1$ (i.e., for bistellar 1-moves) in general.

\begin{eg}$\!\!\!${\rm {\bf :}}  \label{wpm27}
{\rm Some weak $2$-pseudomanifolds on at most seven vertices. The
degree sequences are presented parenthetically below the figures.
Each of $S_1, \dots, S_9$ triangulates the 2-sphere, each of
$R_1, \dots, R_4$ triangulates the real projective plane and $T$
triangulates the torus. Observe that $P_1$, $P_2$ are not
pseudomanifolds. }
\end{eg}


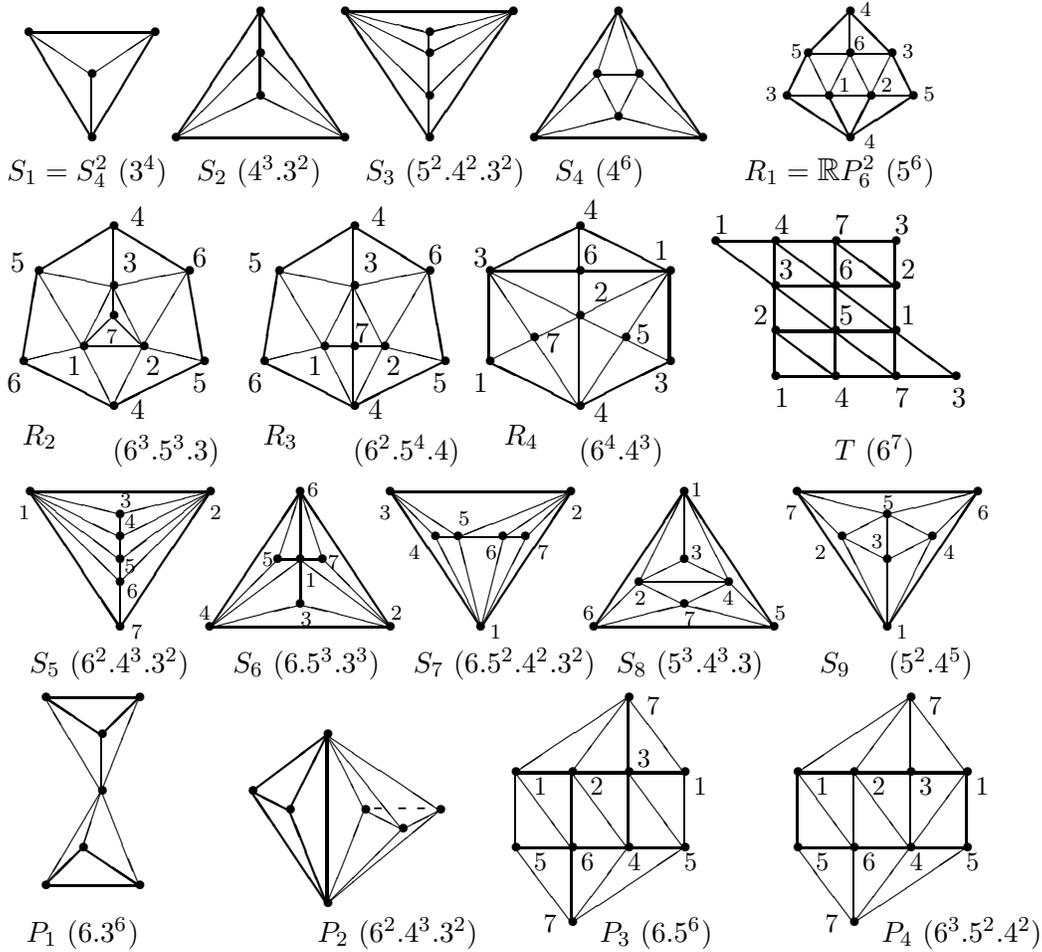
\setlength{\unitlength}{2.8mm}
\begin{picture}(45,9.5)(0,-1)


\put(0.77,7){$_{\bullet}$}
\put(6.77,7){$_{\bullet}$}    
\put(3.77,5){$_{\bullet}$} \put(3.77,2){$_{\bullet}$}

\thicklines

\put(1,7){\line(3,-5){3}} \put(1,7){\line(1,0){6}}
\put(4,2){\line(3,5){3}}

\thinlines

\put(4,2){\line(0,1){3}} \put(4,5){\line(-3,2){3}}
\put(4,5){\line(3,2){3}}

\put(0,0){\mbox{$S_{1}=S^2_{4}$ $(3^4)$}}


\put(7.77,2){$_{\bullet}$} \put(15.77,2){$_{\bullet}$}
\put(11.77,4){$_{\bullet}$} \put(11.77,6){$_{\bullet}$}
\put(11.77,8){$_{\bullet}$}

\thicklines

\put(8,2){\line(2,3){4}} \put(16,2){\line(-2,3){4}}
\put(8,2){\line(1,0){8}}

\thinlines

\put(8,2){\line(2,1){4}} \put(8,2){\line(1,1){4}}

\put(16,2){\line(-2,1){4}} \put(16,2){\line(-1,1){4}}

\put(12,4){\line(0,1){4}}

\put(9,0){\mbox{$S_{2}$ $(4^3.3^2)$}}


\put(19.8,2){$_{\bullet}$} \put(19.8,4){$_{\bullet}$}
\put(19.8,6){$_{\bullet}$} \put(19.8,7){$_{\bullet}$}
\put(15.79,8){$_{\bullet}$} \put(23.77,8){$_{\bullet}$}

\thicklines

\put(20,2){\line(2,3){4}} \put(20,2){\line(-2,3){4}}
\put(16,8){\line(1,0){8}}

\thinlines

\put(20,2){\line(0,1){5}} \put(16,8){\line(4,-1){4}}
\put(16,8){\line(2,-1){4}} \put(16,8){\line(1,-1){4}}
\put(24,8){\line(-4,-1){4}} \put(24,8){\line(-2,-1){4}}
\put(24,8){\line(-1,-1){4}}

\put(17,0){\mbox{$S_{3}$ $(5^2.4^2.3^2)$}}


\put(24.77,2){$_{\bullet}$} \put(32.77,2){$_{\bullet}$}
\put(27.77,5){$_{\bullet}$} \put(29.77,5){$_{\bullet}$}
\put(28.77,8){$_{\bullet}$} \put(28.77,3){$_{\bullet}$}

\thicklines

\put(25,2){\line(2,3){4}} \put(33,2){\line(-2,3){4}}
\put(25,2){\line(1,0){8}}

\thinlines

\put(25,2){\line(1,1){3}} \put(25,2){\line(4,1){4}}
\put(33,2){\line(-4,1){4}} \put(29,3){\line(-1,2){1}}
\put(29,3){\line(1,2){1}} \put(33,2){\line(-1,1){3}}
\put(28,5){\line(1,0){2}} \put(28,5){\line(1,3){1}}
\put(30,5){\line(-1,3){1}}

\put(26,0){\mbox{$S_{4}$ $(4^6)$}}


\put(39.77,2){$_{\bullet}$}
\put(36.77,4){$_{\bullet}$}    
\put(38.77,4){$_{\bullet}$} \put(40.77,4){$_{\bullet}$}
\put(42.77,4){$_{\bullet}$} \put(37.77,6){$_{\bullet}$}
\put(39.77,6){$_{\bullet}$} \put(41.77,6){$_{\bullet}$}
\put(39.77,8){$_{\bullet}$}

\put(40.7,1.7){$_{4}$} \put(36,4){$_{3}$} \put(39.5,4.4){$_{1}$}
\put(41.5,4.4){$_{2}$} \put(43.5,4){$_{5}$}

\put(37.2,6){$_{5}$} \put(40.2,6.4){$_{6}$} \put(42.5,6){$_{3}$}
\put(40.5,8){$_{4}$}

\thicklines

\put(40,2){\line(-3,2){3}} \put(40,2){\line(-1,2){1}}
\put(40,2){\line(3,2){3}} \put(40,2){\line(1,2){1}}
\put(37,4){\line(1,2){1}} \put(43,4){\line(-1,2){1}}
\put(38,6){\line(1,1){2}} \put(42,6){\line(-1,1){2}}

\thinlines

\put(37,4){\line(1,0){6}} \put(38,6){\line(1,0){4}}
\put(39,4){\line(-1,2){1}} \put(39,4){\line(1,2){1}}
\put(41,4){\line(-1,2){1}} \put(41,4){\line(1,2){1}}
\put(40,6){\line(0,1){2}}

\put(35,0){\mbox{$R_{1}=\RR P^2_6$ $(5^6)$}}

\end{picture}



\setlength{\unitlength}{2mm}
\begin{picture}(62,17.5)(1,-4)


\put(7.75,0){$_{\bullet}$} \put(1.75,3){$_{\bullet}$}
\put(13.75,3){$_{\bullet}$}
\put(7.75,6){$_{\bullet}$}    
\put(5.75,4){$_{\bullet}$} \put(9.75,4){$_{\bullet}$}
\put(7.75,8){$_{\bullet}$} \put(2.75,9){$_{\bullet}$}
\put(12.75,9){$_{\bullet}$} \put(7.75,12){$_{\bullet}$}

\put(9.2,-0.8){$ {4}$} \put(1,1){$ {6}$} \put(1.2,9){$ {5}$}
\put(9.2,12){$ {4}$} \put(13.3,9.3){$ {6}$} \put(13.3,1){$ {5}$}
\put(5,2.1){$ {1}$} \put(10.2,2.1){$ {2}$} \put(8.7,9){$ {3}$}
\put(7.6,4.5){$_{7}$}

\thicklines

\put(8,0){\line(-2,1){6}} \put(2,3){\line(1,6){1}}
\put(3,9){\line(5,3){5}} \put(8,12){\line(5,-3){5}}
\put(13,9){\line(1,-6){1}} \put(14,3){\line(-2,-1){6}}

\thinlines

\put(6,4){\line(-4,-1){4}} \put(6,4){\line(-3,5){3}}
\put(6,4){\line(1,2){2}} \put(6,4){\line(1,0){4}}
\put(6,4){\line(1,-2){2}} \put(10,4){\line(-1,-2){2}}
\put(10,4){\line(-1,2){2}} \put(10,4){\line(3,5){3}}
\put(10,4){\line(4,-1){4}} \put(8,8){\line(-5,1){5}}
\put(8,8){\line(0,1){4}} \put(8,8){\line(5,1){5}}
\put(6,4){\line(1,1){2}} \put(10,4){\line(-1,1){2}}
\put(8,8){\line(0,-1){2}}

\put(2,-2.5){\mbox{${R}_{2}$}} \put(8,-3.5){$(6^3.5^3.3)$}


\put(23.75,0){$_{\bullet}$} \put(17.75,3){$_{\bullet}$}
\put(29.75,3){$_{\bullet}$}
\put(23.75,4){$_{\bullet}$}    
\put(21.75,4){$_{\bullet}$} \put(25.75,4){$_{\bullet}$}
\put(23.75,8){$_{\bullet}$} \put(18.75,9){$_{\bullet}$}
\put(28.75,9){$_{\bullet}$} \put(23.75,12){$_{\bullet}$}

\put(25,-0.8){$ {4}$} \put(17,1){$ {6}$} \put(17,9){$ {5}$}
\put(25.2,12){$ {4}$} \put(29,9.5){$ {6}$} \put(29.3,1){$ {5}$}
\put(21,2.1){$ {1}$} \put(26.2,2.1){$ {2}$} \put(24.7,9){$ {3}$}
\put(24.1,4.3){$ {7}$}

\thicklines

\put(24,0){\line(-2,1){6}} \put(18,3){\line(1,6){1}}
\put(19,9){\line(5,3){5}} \put(24,12){\line(5,-3){5}}
\put(29,9){\line(1,-6){1}} \put(30,3){\line(-2,-1){6}}

\thinlines

\put(22,4){\line(-4,-1){4}} \put(22,4){\line(-3,5){3}}
\put(22,4){\line(1,2){2}} \put(22,4){\line(1,-2){2}}
\put(26,4){\line(-1,-2){2}} \put(26,4){\line(-1,2){2}}
\put(26,4){\line(3,5){3}} \put(26,4){\line(4,-1){4}}
\put(24,8){\line(-5,1){5}} \put(24,8){\line(0,1){4}}
\put(24,8){\line(5,1){5}} \put(24,8){\line(0,-1){4}}
\put(24,0){\line(0,1){4}} \put(22,4){\line(1,0){4}}

\put(18,-2.5){\mbox{${R}_{3}$}} \put(24,-3.5){$(6^2.5^4.4)$}


\put(32.74,3){$_{\bullet}$} \put(32.74,9){$_{\bullet}$}
\put(38.74,0){$_{\bullet}$}    
\put(38.74,6){$_{\bullet}$} \put(38.74,9){$_{\bullet}$}
\put(38.74,12){$_{\bullet}$} \put(44.74,3){$_{\bullet}$}
\put(44.74,9){$_{\bullet}$}

\put(35.71,4.57){$_{\bullet}$} \put(41.78,4.54){$_{\bullet}$}

\put(32,1){$ {1}$} \put(32,9.3){$ {3}$} \put(40,-1){$ {4}$}
\put(40,7){$ {2}$} \put(39.3,9.5){$ {6}$} \put(39.3,12.5){$ {4}$}
\put(36.8,3.6){$ {7}$} \put(42.8,4.2){$ {5}$} \put(44,1){$ {3}$}
\put(44,9.8){$ {1}$}

\thicklines

\put(33,3){\line(0,1){6}} \put(39,0){\line(-2,1){6}}
\put(39,0){\line(2,1){6}} \put(33,9){\line(2,1){6}}
\put(45,9){\line(0,-1){6}} \put(45,9){\line(-2,1){6}}

\thinlines

\put(33,3){\line(2,1){12}} \put(39,0){\line(-2,3){6}}
\put(39,0){\line(0,1){12}} \put(39,0){\line(2,3){6}}
\put(33,9){\line(2,-1){12}} \put(33,9){\line(1,0){12}}

\put(34,-2.5){\mbox{$R_{4}$}} \put(39,-3.5){$(6^4.4^3)$}


\put(51.73,2){$_{\bullet}$} \put(55.73,2){$_{\bullet}$}
\put(59.73,2){$_{\bullet}$} \put(63.73,2){$_{\bullet}$}

\put(59.73,5){$_{\bullet}$} \put(51.73,5){$_{\bullet}$}
\put(55.73,5){$_{\bullet}$}

\put(59.73,8){$_{\bullet}$} \put(51.73,8){$_{\bullet}$}
\put(55.73,8){$_{\bullet}$}

\put(59.73,11){$_{\bullet}$} \put(47.73,11){$_{\bullet}$}
\put(51.73,11){$_{\bullet}$} \put(55.73,11){$_{\bullet}$}

\put(52,0){$ {1}$} \put(56,0){$ {4}$} \put(60,0){$ 7$}
\put(63.8,0){$ {3}$}

\put(50.6,5.5){$ {2}$} \put(56.3,5.5){$ {5}$} \put(60.4,5.5){$
{1}$}

\put(52.3,8.5){$ {3}$} \put(56.3,8.5){$ {6}$} \put(60.4,8.5){$
{2}$}

\put(48,11.5){$ {1}$} \put(52,11.5){$ {4}$} \put(56,11.5){$ 7$}
\put(60,11.5){$ {3}$}

\thicklines

\put(52,2){\line(1,0){12}} \put(52,5){\line(1,0){8}}
\put(52,8){\line(1,0){8}} \put(48,11){\line(1,0){12}}

\put(52,2){\line(0,1){9}} \put(56,2){\line(0,1){9}}
\put(60,2){\line(0,1){9}}

\put(52,5){\line(4,-3){4}} \put(48,11){\line(4,-3){12}}
\put(52,11){\line(4,-3){12}} \put(56,11){\line(4,-3){4}}

\put(56,-3.5){\mbox{${T}$ $(6^7)$}}

\end{picture}



\setlength{\unitlength}{3mm}
\begin{picture}(43,10)(-1,-1)


\put(-0.23,8){$_{\bullet}$}
\put(7.77,8){$_{\bullet}$}    
\put(3.77,7){$_{\bullet}$} \put(3.77,6){$_{\bullet}$}
\put(3.77,5){$_{\bullet}$} \put(3.77,4){$_{\bullet}$}
\put(3.77,2){$_{\bullet}$}

\put(-0.5,7){$_{1}$} \put(8,7){$_{2}$} \put(4,7.5){$_{3}$}
\put(4.5,1.7){$_{7}$} \put(4.3,3.5){$_6$} \put(4.2,4.6){$_5$}
\put(4.2,6.6){$_4$}

\thicklines

\put(0,8){\line(2,-3){4}} \put(0,8){\line(1,0){8}}
\put(4,2){\line(2,3){4}}

\thinlines

\put(4,2){\line(0,1){5}} \put(0,8){\line(4,-1){4}}
\put(0,8){\line(2,-1){4}} \put(0,8){\line(4,-3){4}}
\put(0,8){\line(1,-1){4}} \put(8,8){\line(-4,-1){4}}
\put(8,8){\line(-2,-1){4}} \put(8,8){\line(-4,-3){4}}
\put(8,8){\line(-1,-1){4}}

\put(0,0){\mbox{$S_{5}$ $(6^2.4^3.3^2)$}}


\put(7.73,2){$_{\bullet}$} \put(15.73,2){$_{\bullet}$}
\put(11.73,3){$_{\bullet}$} \put(11.73,5){$_{\bullet}$}
\put(10.73,5){$_{\bullet}$} \put(12.73,5){$_{\bullet}$}
\put(11.73,8){$_{\bullet}$}

\put(7.5,2.5){$_{4}$} \put(16,2.5){$_{2}$} \put(12,2.25){$_{3}$}
\put(12.2,4){$_{1}$} \put(12.3,8){$_{6}$} \put(10.3,4.9){$_{5}$}
\put(13.2,4.8){$_{7}$}

\thicklines

\put(8,2){\line(2,3){4}} \put(16,2){\line(-2,3){4}}
\put(8,2){\line(1,0){8}}

\thinlines

\put(8,2){\line(4,1){4}} \put(8,2){\line(4,3){4}}
\put(8,2){\line(1,1){3}} \put(16,2){\line(-4,1){4}}
\put(16,2){\line(-4,3){4}} \put(16,2){\line(-1,1){3}}
\put(12,3){\line(0,1){5}}

\put(11,5){\line(1,0){2}} \put(11,5){\line(1,3){1}}
\put(13,5){\line(-1,3){1}}

\put(9,0){\mbox{$S_{6}$ $(6.5^3.3^3)$}}


\put(17.73,6){$_{\bullet}$} \put(18.73,6){$_{\bullet}$}
\put(20.73,6){$_{\bullet}$} \put(21.73,6){$_{\bullet}$}
\put(15.73,8){$_{\bullet}$} \put(23.73,8){$_{\bullet}$}
\put(19.73,2){$_{\bullet}$}

\put(15.5,7){$_{3}$} \put(24,7){$_{2}$} \put(20.4,1.7){$_{1}$}
\put(16.8,5.2){$_{4}$} \put(19,6.7){$_{5}$}
\put(20.2,5.2){$_{6}$} \put(22.5,5.2){$_{7}$}

\thicklines

\put(20,2){\line(2,3){4}} \put(20,2){\line(-2,3){4}}
\put(16,8){\line(1,0){8}}

\thinlines

\put(20,2){\line(-1,2){2}} \put(20,2){\line(-1,4){1}}
\put(20,2){\line(1,2){2}} \put(20,2){\line(1,4){1}}

\put(24,8){\line(-5,-2){5}} \put(16,8){\line(3,-2){3}}
\put(16,8){\line(1,-1){2}} \put(24,8){\line(-3,-2){3}}
\put(24,8){\line(-1,-1){2}} \put(18,6){\line(1,0){4}}

\put(17,0){\mbox{$S_{7}$ $(6.5^2.4^2.3^2)$}}


\put(24.75,2){$_{\bullet}$} \put(32.75,2){$_{\bullet}$}
\put(28.75,3){$_{\bullet}$} \put(26.75,4){$_{\bullet}$}
\put(30.73,4){$_{\bullet}$} \put(28.76,8){$_{\bullet}$}
\put(28.76,5){$_{\bullet}$}

\put(24.5,2.5){$_{6}$} \put(33,2.5){$_{5}$} \put(29,2.25){$_{7}$}
\put(26.8,3.2){$_{2}$} \put(29.3,8){$_1$} \put(29.3,5.2){$_3$}
\put(30.7,3.2){$_4$}

\thicklines

\put(25,2){\line(2,3){4}} \put(33,2){\line(-2,3){4}}
\put(25,2){\line(1,0){8}}

\thinlines

\put(25,2){\line(1,1){2}} \put(25,2){\line(4,1){4}}
\put(33,2){\line(-4,1){4}} \put(33,2){\line(-1,1){2}}
\put(29,3){\line(2,1){2}} \put(29,3){\line(-2,1){2}}
\put(27,4){\line(1,0){4}} \put(27,4){\line(2,1){2}}
\put(27,4){\line(1,2){2}} \put(31,4){\line(-2,1){2}}
\put(31,4){\line(-1,2){2}} \put(29,5){\line(0,1){3}}

\put(26,0){\mbox{$S_{8}$ $(5^3.4^3.3)$}}


\put(37.76,2){$_{\bullet}$} \put(37.76,5){$_{\bullet}$}
\put(37.76,7){$_{\bullet}$} \put(35.76,6){$_{\bullet}$}
\put(39.76,6){$_{\bullet}$} \put(41.76,8){$_{\bullet}$}
\put(33.76,8){$_{\bullet}$}

\put(33.5,7){$_{7}$} \put(42,7){$_{6}$} \put(38.4,1.7){$_{1}$}
\put(34.8,5.2){$_{2}$} \put(37.3,5.7){$_3$} \put(37.8,7.4){$_5$}
\put(40.5,5.2){$_4$}

\thicklines

\put(38,2){\line(2,3){4}} \put(38,2){\line(-2,3){4}}
\put(34,8){\line(1,0){8}}

\thinlines

\put(38,2){\line(0,1){5}} \put(38,2){\line(1,2){2}}
\put(38,2){\line(-1,2){2}} \put(38,5){\line(-2,1){2}}
\put(38,5){\line(2,1){2}} \put(38,7){\line(4,1){4}}
\put(38,7){\line(2,-1){2}} \put(38,7){\line(-4,1){4}}
\put(38,7){\line(-2,-1){2}} \put(36,6){\line(-1,1){2}}
\put(40,6){\line(1,1){2}}

\put(35,0){\mbox{$S_9$ ~~ $(5^2.4^5)$}}

\end{picture}



\setlength{\unitlength}{2.5mm}
\begin{picture}(46,14)(-12,0)


\put(-7.23,12){$_{\bullet}$} \put(-10.23,14){$_{\bullet}$}
\put(-5.23,14){$_{\bullet}$} \put(-7.23,9){$_{\bullet}$}
\put(-8.23,6){$_{\bullet}$} \put(-5.23,4){$_{\bullet}$}
\put(-10.23,4){$_{\bullet}$}

\thicklines

\put(-7,12){\line(-3,2){3}} \put(-7,12){\line(1,1){2}}
\put(-10,14){\line(1,0){5}}

\put(-8,6){\line(3,-2){3}} \put(-10,4){\line(1,1){2}}
\put(-10,4){\line(1,0){5}}

\thinlines

\put(-7,9){\line(-3,5){3}} \put(-7,9){\line(2,5){2}}
\put(-7,9){\line(0,1){3}} \put(-7,9){\line(2,-5){2}}
\put(-7,9){\line(-1,-3){1}} \put(-7,9){\line(-3,-5){3}}

\put(-11,1){\mbox{$P_1$ $ (6.3^6)$}}


\put(4.77,3){$_{\bullet}$} \put(8.77,7){$_{\bullet}$}
\put(2.77,8){$_{\bullet}$} \put(6.77,8){$_{\bullet}$}
\put(10.77,8){$_{\bullet}$} \put(0.77,9){$_{\bullet}$}
\put(4.77,12){$_{\bullet}$}

\put(7.5,7.78){-\, - - -} 

\thicklines

\put(1,9){\line(4,3){4}} \put(3,8){\line(1,2){2}}
\put(3,8){\line(-2,1){2}} \put(5,3){\line(-2,3){4}}
\put(5,3){\line(-2,5){2}} \put(5,3){\line(0,1){9}}

\thinlines

\put(9,7){\line(2,1){2}} \put(9,7){\line(-2,1){2}}

\put(5,3){\line(2,5){2}} \put(5,3){\line(1,1){4}}
\put(5,3){\line(6,5){6}} \put(5,12){\line(1,-2){2}}
\put(5,12){\line(4,-5){4}} \put(5,12){\line(3,-2){6}}

\put(4.5,1){\mbox{$P_2$ $(6^2.4^3.3^2)$}}


\put(16,8.8){\small 1} \put(19,8.8){\small 2}
\put(21.5,10.3){\small 3} \put(24.5,8.8){\small 1}
\put(16,4.8){\small 5} \put(18.5,4.8){\small 6}
\put(21,4.8){\small 4} \put(24,4.8){\small 5}
\put(16.5,1.5){\small 7} \put(22,13.1){\small 7}

\put(17.77,2){$_{\bullet}$} \put(14.77,6){$_{\bullet}$}
\put(17.77,6){$_{\bullet}$} \put(20.77,6){$_{\bullet}$}
\put(23.77,6){$_{\bullet}$} \put(14.77,10){$_{\bullet}$}
\put(17.77,10){$_{\bullet}$} \put(20.77,10){$_{\bullet}$}
\put(23.77,10){$_{\bullet}$} \put(20.77,14){$_{\bullet}$}

\thicklines

\put(15,6){\line(1,0){9}} \put(15,10){\line(1,0){9}}

\thinlines

\put(18,2){\line(-3,4){3}} \put(18,2){\line(0,1){4}}
\put(18,2){\line(3,4){3}} \put(18,2){\line(3,2){6}}

\put(21,14){\line(-3,-2){6}} \put(21,14){\line(-3,-4){3}}
\put(21,14){\line(0,-1){4}} \put(21,14){\line(3,-4){3}}

\put(15,10){\line(3,-4){3}} \put(18,10){\line(3,-4){3}}
\put(21,10){\line(3,-4){3}}

\put(15,10){\line(0,-1){4}} \put(18,10){\line(0,-1){4}}
\put(21,10){\line(0,-1){4}} \put(24,10){\line(0,-1){4}}

\put(19.5,1){\mbox{$P_3$ $(6.5^6)$}}


\put(31,8.8){\small 1} \put(34,8.8){\small 2}
\put(36.5,8.8){\small 3} \put(39.5,8.8){\small 1}
\put(31,4.8){\small 5} \put(33.5,4.8){\small 6}
\put(36,4.8){\small 4} \put(39,4.8){\small 5}
\put(31.5,1.5){\small 7} \put(37,13.1){\small 7}

\put(32.77,2){$_{\bullet}$} \put(29.77,6){$_{\bullet}$}
\put(32.77,6){$_{\bullet}$} \put(35.77,6){$_{\bullet}$}
\put(38.77,6){$_{\bullet}$} \put(29.77,10){$_{\bullet}$}
\put(32.77,10){$_{\bullet}$} \put(35.77,10){$_{\bullet}$}
\put(38.77,10){$_{\bullet}$} \put(35.77,14){$_{\bullet}$}

\thicklines

\put(30,6){\line(1,0){9}} \put(30,10){\line(1,0){9}}

\thinlines

\put(33,2){\line(-3,4){3}} \put(33,2){\line(0,1){4}}
\put(33,2){\line(3,4){3}} \put(33,2){\line(3,2){6}}

\put(36,14){\line(-3,-2){6}} \put(36,14){\line(-3,-4){3}}
\put(36,14){\line(0,-1){4}} \put(36,14){\line(3,-4){3}}

\put(30,10){\line(3,-4){3}} \put(33,10){\line(3,-4){3}}
\put(36,6){\line(3,4){3}}

\put(30,10){\line(0,-1){4}} \put(33,10){\line(0,-1){4}}
\put(36,10){\line(0,-1){4}} \put(39,10){\line(0,-1){4}}

\put(34.5,1){\mbox{$P_4$ $(6^3.5^2.4^2)$}}

\end{picture}


\medskip

We know that if $K$ is a weak 2-pseudomanifold with at most six
vertices then $K$ is isomorphic to $S_1, \dots, S_4$ or $R_1$
(cf. \cite{bd2}). In \cite{d}, we have seen the following.

\begin{prop}$\!\!\!${\bf .} \label{w2mfd}
There are exactly $13$ distinct $2$-dimensional weak
pseudomanifolds on $7$ vertices, namely, $S_{5}, \dots, S_{9}$,
$R_{2}, \dots, R_{4}$, $T$, $P_1, \dots, P_3$ and $P_4$.
\end{prop}

\section{Examples}

We identify a weak pseudomanifold with the set of facets in it.

\begin{eg}$\!\!\!${\rm {\bf :}}  \label{m38n}
{\rm These four neighbourly 8-vertex combinatorial 3-manifolds
were found by Gr\"{u}nbaum and Sreedharan (in \cite{gs}, these
are denoted by $P^8_{35}$, $P^8_{36}$, $P^8_{37}$ and ${\cal M}$
respectively). It follows from Lemma \ref{le3.2} that these are
combinatorial 3-spheres. It was shown in \cite{gs} that the first
three of these are polytopal 3-spheres and the last one is a
non-polytopal sphere.}
\begin{eqnarray*}
S^{\, 3}_{8,35} & = & \{1234, 1267, 1256, 1245, 2345, 2356, 2367,
3467,
3456, 4567, \\
&& \quad 1238, 1278, 2378, 1348, 3478, 1458, 4578, 1568, 1678,
5678\}, \\
S^{\, 3}_{8,36} & = & \{1234, 1256, 1245, 1567, 2345, 2356, 2367,
3467,
3456, 4567, \\
&& \quad 1268, 1678, 2678, 1238, 2378, 1348, 3478, 1458, 1578,
4578\}, \\
S^{\, 3}_{8,37} & = & \{1234, 1256, 1245, 1457, 2345, 2356, 2367,
3467,
3456, 4567, \\
&& \quad 1568, 1578, 5678, 1268, 2678, 1238, 2378, 1348, 1478,
3478\}, \\
S^{\, 3}_{8,38} & = & \{1234, 1237, 1267, 1347, 1567, 2345, 2367,
3467,
3456, 4567, \\
&& \quad 2358, 2368, 3568, 1268, 1568, 1248, 2458, 1478, 1578,
4578\}.
\end{eqnarray*}
\end{eg}

\begin{lemma}$\!\!\!${\bf .} \label{le3.1}
$S^3_{8,i} \not \cong S^3_{8,j}$ for $35 \leq i<j\leq 38$.
\end{lemma}

\noindent {\bf Proof.} Observe that $G_6(S^3_{8,35}) = C_8(1, 2,
\dots, 8)$, $G_6(S^3_{8,36}) = (V, \{23, 34, 45, 67, 78, 81\})$,
$G_6(S^3_{8,37})$ $ = (V, \{23, 34, 56, 78, 81\})$ and
$G_6(S^3_{8,38}) = (V, \{17, 23, 58\})$, where $V = \{1, \dots,
8\}$. Since $K \cong L$ implies $G_6(K) \cong G_6(L)$, $S^3_{8,i}
\not\cong S^3_{8,j}$, for $35 \leq i< j \leq 38$. \hfill $\Box$

\begin{eg}$\!\!\!${\rm {\bf :}}  \label{m38nn}
{\rm Some non-neighbourly 8-vertex combinatorial 3-manifolds. It
follows from Lemma \ref{le3.2} that these are combinatorial
3-spheres. For $1 \leq i \leq 34$, the sphere $S^3_{8, i}$ is
isomorphic to the polytopal sphere $P^{\,8}_i$ in \cite{gs} and
the sphere $S^3_{8, 39}$ is isomorphic to the non-polytopal
sphere found by Barnette in \cite{b1}. We consecutively define\,:}
\begin{eqnarray*}
  S^{\, 3}_{8,39} = \kappa_{46}(S^{\, 3}_{8, 38}), &&
  S^{\, 3}_{8, 33} = \kappa_{27}(S^{\, 3}_{8, 37}), \hspace{4mm}
  S^{\, 3}_{8, 32} = \kappa_{48}(S^{\, 3}_{8, 37}), \hspace{4mm}
  S^{\, 3}_{8, 31} = \kappa_{58}(S^{\, 3}_{8, 37}), \\
 S^{\, 3}_{8, 30} = \kappa_{24}(S^{\, 3}_{8, 37}), &&
 S^{\, 3}_{8, 29} = \kappa_{27}(S^{\, 3}_{8, 31}), \hspace{4mm}
 S^{\, 3}_{8, 28} = \kappa_{24}(S^{\, 3}_{8, 31}), \hspace{4mm}
 S^{\, 3}_{8, 27} = \kappa_{13}(S^{\, 3}_{8, 31}), \\
 S^{\, 3}_{8, 25} = \kappa_{57}(S^{\, 3}_{8, 31}), &&
 S^{\, 3}_{8, 24} = \kappa_{48}(S^{\, 3}_{8, 31}), \hspace{4mm}
 S^{\, 3}_{8, 23} = \kappa_{35}(S^{\, 3}_{8, 31}),\\
 S^{\, 3}_{8, 26} = \kappa_{46}(S^{\, 3}_{8, 27}), &&
 S^{\, 3}_{8, 22} = \kappa_{24}(S^{\, 3}_{8, 25}), \hspace{4mm}
 S^{\, 3}_{8, 21} = \kappa_{68}(S^{\, 3}_{8, 25}), \hspace{4mm}
 S^{\, 3}_{8, 20} = \kappa_{48}(S^{\, 3}_{8, 25}), \\
 S^{\, 3}_{8, 19} =  \kappa_{17}(S^{\, 3}_{8, 25}), &&
 S^{\, 3}_{8, 18} = \kappa_{27}(S^{\, 3}_{8, 25}), \hspace{4mm}
 S^{\, 3}_{8, 12} = \kappa_{15}(S^{\, 3}_{8, 25}), \hspace{4mm}
 S^{\, 3}_{8, 11} = \kappa_{35}(S^{\, 3}_{8, 25}), \\
 S^{\, 3}_{8, 17} = \kappa_{24}(S^{\, 3}_{8, 19}), &&
 S^{\, 3}_{8,34} = \kappa_{27}(S^{\, 3}_{8, 26})
 \hspace{0.5mm} = \hspace{0.5mm} S^{\, 0}_3(1, 3) \ast S^{\, 0}_3(2, 7)
 \ast S^{\, 0}_3(4, 6) \ast S^{\, 0}_3(5, 8),\\
 S^{\, 3}_{8, 16} = \kappa_{13}(S^{\, 3}_{8, 19}), &&
 S^{\, 3}_{8, 15} = \kappa_{28}(S^{\, 3}_{8, 18}), \hspace{4mm}
 S^{\, 3}_{8, 14} =  \kappa_{47}(S^{\, 3}_{8, 20}), \hspace{4mm}
 S^{\, 3}_{8, 10} = \kappa_{15}(S^{\, 3}_{8, 19}), \\
 S^{\, 3}_{8, 9} = \kappa_{35}(S^{\, 3}_{8, 19}), && \hspace{1.5mm}
 S^{\, 3}_{8, 8} = \kappa_{47}(S^{\, 3}_{8, 19}), \hspace{4mm}
 S^{\, 3}_{8, 13} = \kappa_{38}(S^{\, 3}_{8, 16}), \hspace{5.5mm}
 S^{\, 3}_{8, 7} = \kappa_{24}(S^{\, 3}_{8, 8}), \\
 S^{\, 3}_{8, 6} = \kappa_{35}(S^{\, 3}_{8, 8}), \hspace{1.2mm} && \hspace{1.5mm}
 S^{\, 3}_{8, 5} = \kappa_{48}(S^{\, 3}_{8, 8}), \hspace{6.5mm}
 S^{\, 3}_{8, 4} = \kappa_{15}(S^{\, 3}_{8, 8}), \\
 S^{\, 3}_{8, 3} = \kappa_{48}(S^{\, 3}_{8, 4}), \hspace{1.5mm} && \hspace{1.5mm}
 S^{\, 3}_{8, 2} = \kappa_{48}(S^{\, 3}_{8, 6}), \hspace{7mm}
 S^{\, 3}_{8, 1} = \kappa_{16}(S^{\, 3}_{8, 4}).
\end{eqnarray*}
\end{eg}


\begin{lemma}$\!\!\!${\bf .} \label{le3.2} $(a)$ $S^3_{8,i}
\approx S^3_{8,j}$, for $1 \leq i, j \leq 39$, $(b)$ $S^3_{8, m}$
is a combinatorial $3$-sphere, for $1 \leq m \leq 39$ and $(c)$
$S^3_{8,k} \not \cong S^3_{8,l}$, for $1 \leq k < l \leq 39$.
\end{lemma}

\noindent {\bf Proof.} For $0\leq i\leq 6$, let ${\cal S}_i$
denote the set of $S^3_{8,j}$'s with $i$ non-edges. Then ${\cal
S}_0 = \{S^3_{8,35}, S^3_{8,36}, S^3_{8,37}$, $S^3_{8,38}\}$,
${\cal S}_1 = \{S^3_{8,30}, S^3_{8,31}, S^3_{8,32}, S^3_{8,33},
S^3_{8,39}\}$, ${\cal S}_2 = \{S^3_{8, 23}, S^3_{8, 24}, S^3_{8,
25}, S^3_{8, 27}$, $S^3_{8, 28}, S^3_{8, 29}\}$, ${\cal S}_3 =
\{S^3_{8,11}$, $S^3_{8,12}, S^3_{8,18}, S^3_{8,19}, S^3_{8,20},
S^3_{8,21}, S^3_{8,22}, S^3_{8,26}\}$, ${\cal S}_4 = \{S^3_{8,8},
S^3_{8,9}, S^3_{8,10}$, $S^3_{8,14}, S^3_{8,15}, S^3_{8,16},
S^3_{8,17}, S^3_{8,34}\}$, ${\cal S}_5 $ $ = \{S^3_{8,4},
S^3_{8,5}, S^3_{8,6}, S^3_{8,7}, S^3_{8,13}\}$ and ${\cal S}_6 =
\{S^3_{8,1}, S^3_{8,2}, S^3_{8,3}\}$.

From the proof of Lemma \ref{S0toS1}, $S^3_{8, 35} \approx S^3_{8,
30} \approx S^3_{8, 36} \approx S^3_{8, 30} \approx S^3_{8, 37}
\approx S^3_{8, 32} \approx S^3_{8, 38}$. Thus, $ S^3_{8, i}
\approx S^3_{8, j}$ for $35 \leq i, j\leq 38$. Now, if $S^3_{8,
i} \in {\cal S}_2 \cup {\cal S}_3 \cup {\cal S}_4 \cup {\cal S}_5
\cup {\cal S}_6$ then, from the definition of $S^3_{8, i}$,
$S^3_{8, i} \approx S^3_{8, 31} \approx S^3_{8, 37}$. This proves
Part $(a)$.

Since $S^3_{8, 34}$ is a join of spheres, $S^3_{8, 34}$ is a
combinatorial 3-sphere. Clearly, if $M \approx N$ and $M$ is a
combinatorial $3$-sphere then $N$ is so. Part $(b)$ now follows
from Part $(a)$.

Since the non-edge graphs of the members of ${\cal S}_6$
(respectively, ${\cal S}_5$) are pairwise non-isomorphic, the
members of ${\cal S}_6$ (respectively, ${\cal S}_5$) are pairwise
non-isomorphic.

For $S^3_{8,i}, S^3_{8,j} \in {\cal S}_4$ ($i < j$) and ${\rm
NEG}(S^3_{8,i}) \cong {\rm NEG}(S^3_{8,j})$ imply $(i, j) = (8,
9)$ or $(14, 15)$. Since $M \cong N$ implies $G_6(M) \cong
G_6(N)$ and $G_6(S^3_{8, 8}) \not \cong G_6(S^3_{8, 9})$,
$G_6(S^3_{8, 14}) \not\cong G_6(S^3_{8, 15})$, the members of
${\cal S}_4$ are pairwise non-isomorphic.

For $S^3_{8,i} \neq S^3_{8,j} \in {\cal S}_3$ and ${\rm
NEG}(S^3_{8,i}) \cong {\rm NEG}(S^3_{8,j})$ imply $\{i, j\} =
\{11, 12\}$ or $18\leq i\neq j \leq 21$. Since the non-edge graph
of a member in $\Sigma_i$ is non-isomorphic to the non-edge graph
of a member of $\Sigma_j$ for $i \neq j$, a member of $\Sigma_i$
is non-isomorphic to a member of $\Sigma_j$. Observe that
$G_6(S^3_{8, 11}) \not \cong G_6(S^3_{8, 12})$ and for $18\leq i<
j\leq 21$, $G_6(S^3_{8, i}) \cong G_6(S^3_{8, j})$ implies $(i,
j) = (18, 19)$. Since $G_3(S^3_{8, 18}) \not \cong G_3(S^3_{8,
19})$ the members of ${\cal S}_3$ are pairwise non-isomorphic.

Since $G_3(S^3_{8, i}) \not \cong G_3(S^3_{8, j})$ for $S^3_{8,i}
\neq S^3_{8,j} \in {\cal S}_2$, the members of ${\cal S}_2$ are
pairwise non-isomorphic. By the same reasoning, the members of
${\cal S}_1$ are pairwise non-isomorphic.

By Lemma \ref{le3.1}, the members of ${\cal S}_0$ are pairwise
non-isomorphic. Since a member of ${\cal S}_i$ is non-isomorphic
to a member of ${\cal S}_j$ for $i\neq j$, the above imply Part
$(c)$. \hfill $\Box$


\begin{eg}$\!\!\!${\rm {\bf :}} \label{pm38}
{\rm Some  8-vertex neighbourly normal 3-pseudomanifolds. }
\begin{eqnarray*}
N_{1} & = & \{1248, 1268, 1348, 1378, 1568, 1578, 2358, 2378,
2458, 2678, 3468, 3568, 4578, 4678, \\
&& 1247, 1257, 1367, 1467, 2347, 2567, 3457, 3567, 1236,
2346, 1345, 1235, 1456, 2456\}, \\
N_{2} & = & \{1248, 2458, 2358, 3568, 3468, 4678, 4578, 1578,
1568, 1268, 2678,  \\
&& \quad 2378, 1378, 1348, 1247, 2457, 2357, 3567, 3467, 1567,
1267, 1347\} = \Sigma_{78}T, \\
N_{3} & = & \{1248, 1268, 1348, 1378, 1568, 1578, 2358, 2378,
2458, 2678, 3468, 3568,  \\
&& \qquad \quad 4578, 4678, 1234, 2347, 2456, 2467, 3456, 3457,
1235, 1256, 1357\}, \\
N_{4} & = & \{1248, 1268, 1348, 1378, 1568, 1578, 2358, 2378,
2458, 2678, 3468, \\
& & \qquad \quad 3568, 4578, 4678, 1245, 1256, 2356, 2367,
3467, 1347, 1457\}, \\
N_{5} & = & \{1258, 1268, 1358, 1378, 1468, 1478, 2368, 2378,
2458, 2478, 3458, 3468,  \\
&& \quad 1257, 1267, 1367, 1457, 2357, 2467, 3457, 3467, 2356,
2456, 1356, 1456\}, \\
N_{6} & = & \{1358, 1378, 1468, 1478, 1568, 2368, 2378, 2458,
2478, 2568, 3458, 3468, \\
&& \quad 1235, 1245, 1457, 1567, 2357, 2567, 3457, 1236, 1246,
1367, 2467, 3467\},\\
N_{7} & = & \{1268, 1258, 1358, 1378, 1478, 1468, 2378, 2368,
2458, 2478, 3468, \\
&& \qquad \quad 3458, 1356, 1367, 2357, 2356, 3467, 3457, 1256,
1467, 2457\}, \\
N_{8} & = & \kappa_{348}(\kappa_{238}(\kappa_{56}(\kappa_{67}
(N_{7})))), ~~ N_{9}  =  \kappa_{235}(\kappa_{67}(N_{7})), \\
N_{10}  & = & \kappa_{148}(\kappa_{67}(N_{7})),  ~~ N_{11}  =
\kappa_{348}(\kappa_{56}((N_{10})), ~~ N_{12}  =
\kappa_{457}(\kappa_{23}((N_{9})), \\ N_{13} & = &
\kappa_{567}(\kappa_{23} ((N_{9})), ~~ N_{14}  =
\kappa_{138}(\kappa_{57}((N_{8}))\cong \Sigma_{78}R_2, ~~ N_{15} =
\kappa_{158}(\kappa_{23}(N_{9})).
\end{eqnarray*}
{\rm All the vertices of $N_{1}$ are singular and their links are
isomorphic to the 7-vertex torus $T$. There are two singular
vertices in $N_{2}$ and their links are isomorphic to $T$. The
singular vertices in $N_{3}$ are  8, 3, 4, 2, 5 and their links
are isomorphic to $T$, $R_2$, $R_2$, $R_3$ and $R_3$ respectively.
There is only one singular vertex in $N_{4}$ whose link is
isomorphic to $T$. All the vertices of $N_{5}$ (respectively,
$N_{6}$) are singular and their links are isomorphic to $R_4$
(respectively, $R_3$). Each of $N_{7}, \dots, N_{15}$ has exactly
two singular vertices and their links  are 7-vertex $\RR P^2$'s.
Thus, each $N_{i}$ is a normal 3-pseudomanifold. }
\end{eg}

It follows from the definition that $N_{i} \approx N_{j}$ for $7
\leq i, j \leq 15$. Here we prove.

\begin{lemma}$\!\!\!${\bf .} \label{le3.3} $(a)$ The geometric carriers
of \, $N_1, N_{2}, N_{3}, N_{4}, N_{5}$ and $N_7$ are distinct
$($non-homeomorphic$)$, $(b)$ $N_{i} \not \approx N_{j}$ for $1
\leq i < j \leq 7$, $(c)$ $N_{5} \sim N_{6}$.
\end{lemma}

\noindent {\bf Proof.} For a normal 3-pseudomanifold $X$, let
$n_s(X)$ denote the number of singular vertices. Clearly, if $M$
and $N$ are two normal 3-pseudomanifolds with homeomorphic
geometric carriers then $(n_s(M), \chi(M)) = (n_s(N)$,
$\chi(N))$. Now, $(n_s(N_{1}), \chi(N_{1})) = (8, 8)$,
$(n_s(N_{2}), \chi(N_{2})) = (2, 2)$, $(n_s(N_{3}), \chi(N_{3}))
= (5, 3)$, $(n_s(N_{4})$, $\chi(N_{4})) = (1, 1)$, $(n_s(N_{5})$,
$\chi(N_{5})) = (8, 4)$, $(n_s(N_{7}), \chi(N_{7})) = (2, 1)$.
This proves $(a)$.

Part $(b)$ follows from the fact that $N_{i}$ is neighbourly and
has no removable edge and hence there is no proper bistellar move
from $N_{i}$, for $1\leq i \leq 6$.

Let $N_{5}^{\prime}$ be obtained from $N_{5}$ by starring a new
vertex 0 in the facet $1358$. Let $N_{5}^{\prime\prime} =
\kappa_{\{0\}} (\kappa_{08} (\kappa_{156} (\kappa_{07}
(\kappa_{03}(\kappa_{035} (\kappa_{68} (\kappa_{02}
(\kappa_{268}(\kappa_{13}(\kappa_{135}(\kappa_{138} (\kappa_{158}
(N_{5}^{\prime})))))))))))))$. Then $N_{5}^{\prime\prime}$ is
isomorphic to $N_{6}$ via the map $(2,3)(5,8)$. This proves
$(c)$. \hfill $\Box$

\begin{lemma}$\!\!\!${\bf .} \label{le3.4}
$N_{k} \not \cong N_{l}$ for $1\leq k<l\leq 15$.
\end{lemma}

\noindent {\bf Proof.} Let $n_s$ be as above. Clearly, if $M$ and
$N$ are two isomorphic 3-pseudomanifolds then $(n_s(M), f_3(M)) =
(n_s(N)$, $f_3(N))$. Now, $(n_s(N_{1}), f_3(N_{1})) = (8, 28)$,
$(n_s(N_{2}), f_3(N_{2})) = (2, 22)$, $(n_s(N_{3}), f_3(N_{3})) =
(5, 23)$,  $(n_s(N_{4}), f_3(N_{4})) = (1, 21)$, $(n_s(N_{5}),
f_3(N_{5})) = (n_s(N_{6})$, $f_3(N_{6})) = (8, 24)$, and
$(n_s(N_{i}), f_3(N_{i})) = (2, 21)$ for $7 \leq i \leq 15$.
Since the links of each vertex in $N_{5}$ is isomorphic to
$R_{4}$ and the links of each vertex in $N_{6}$ is isomorphic to
$R_{3}$, it follows that $N_{5} \not \cong N_{6}$. Thus, $N_{i}
\not \cong N_{j}$ for $1 \leq i \leq 6$, $1 \leq j \leq 15$, $i
\neq j$.

Observe that the singular vertices in $N_i$ are 3 and 8 for $7\leq
i\leq 15$. Moreover, (i) ${\rm lk}_{N_{7}}(3) \cong {\rm
lk}_{N_{7}}(8) \cong R_4$, (ii) ${\rm lk}_{N_{8}}(3) \cong R_4$
and ${\rm lk}_{N_{8}}(8) \cong R_3$, (iii) ${\rm lk}_{N_{9}}(3)
\cong R_2$ and ${\rm lk}_{N_{9}}(8) \cong R_4$, (iv) ${\rm
lk}_{N_{10}}(3) \cong {\rm lk}_{N_{10}}(8) \cong R_3$ and
$\deg_{N_{10}}(38) = 6$, (v) ${\rm lk}_{N_{11}}(3) \cong {\rm
lk}_{N_{11}}(8) \cong R_3$ and $\deg_{N_{11}}(38) = 5$. (vi)
${\rm lk}_{N_{12}}(3) \cong R_2$, ${\rm lk}_{N_{12}}(8) \cong
R_3$ and $G_3(N_{12}) = (V, \{32, 21, 17, 75, 54, 46\})$, (vii)
${\rm lk}_{N_{13}}(3) \cong R_2$, ${\rm lk}_{N_{13}}(8) \cong
R_3$ and $G_3(N_{13}) = (V, \{32, 21$, $17, 75, 56, 67, 64,
42\})$, (viii) ${\rm lk}_{N_{14}}(3) \cong {\rm lk}_{N_{14}}(8)
\cong R_2$ and $\deg_{N_{14}}(38) = 3$. (xi) ${\rm lk}_{N_{15}}(3)
\cong {\rm lk}_{N_{15}}(8) \cong R_2$ and $\deg_{N_{15}}(38) =
6$. These imply that there is no isomorphism between $N_{i}$ and
$N_{j}$ for $7\leq i< j \leq 15$. This completes the proof.
\hfill $\Box$

\begin{eg}$\!\!\!${\rm {\bf :}}  \label{pm38nn}
{\rm Some 8-vertex non-neighbourly normal 3-pseudomanifolds.}
\begin{eqnarray*}
N_{16} = \kappa_{67}(N_{7}), && N_{17} = \kappa_{24}(N_{8}),
\hspace{5mm} N_{18} = \kappa_{238}(\kappa_{56}(\kappa_{67}
(N_{7}))), \hspace{5mm} N_{19} = \kappa_{57}(N_{8}), \\ N_{20} =
\kappa_{56}(N_{10}), && N_{21} = \kappa_{12}(N_{9}), \hspace{5mm}
N_{22} =
\kappa_{14}(N_{11}), \hspace{5mm} N_{23} = \kappa_{23}(N_{9}), \\
N_{24} = \kappa_{38}(N_{14}), && N_{25} = \kappa_{56}(N_{16}),
\hspace{5mm} N_{26} = \kappa_{12}(N_{16}),
\hspace{5mm} N_{27} = \kappa_{56}(N_{17}), \\
N_{28} = \kappa_{57}(N_{18}), && N_{29} = \kappa_{15}(N_{18}),
\hspace{5mm} N_{30} = \kappa_{12}(N_{23}), \hspace{5mm} N_{31} =
\kappa_{24}(N_{22}), \\ N_{32} = \kappa_{24}(N_{26}), && N_{33} =
\kappa_{57}(N_{25}),\hspace{5mm} N_{34} = \kappa_{45}(N_{28}),
\hspace{5mm} N_{35} = \kappa_{58}(N_{29}).
\end{eqnarray*}
\end{eg}

\begin{lemma}$\!\!\!${\bf .} \label{le3.5}
$(a)$ $N_{i} \not \cong N_{j}$ for $1 \leq i < j \leq 35$ and
$(b)$ $N_{k} \approx N_{l}$, for $7 \leq k,  l \leq 35$.
\end{lemma}

\noindent {\bf Proof.} For $0\leq i\leq 3$, let ${\cal N}_i$
denote the set of 3-pseudomanifolds defined in Examples 4 and 5
with $i$ non-edges. Then ${\cal N}_0 = \{N_{1}, \dots, N_{15}\}$,
${\cal N}_1 = \{N_{16}, \dots, N_{24}\}$, ${\cal N}_2 = \{N_{25},
\dots, N_{31}\}$ and ${\cal N}_3 = \{N_{32}, \dots, N_{35}\}$. The
singular vertices in $N_{i}$ are 3 and 8 for $7\leq i\leq 35$.

By Lemma \ref{le3.4}, the members of ${\cal N}_0$ are pairwise
non-isomorphic.

Observe that (i) ${\rm lk}_{N_{16}}(3) \cong R_4$ and ${\rm
lk}_{N_{16}}(8) \cong R_3$, (ii) ${\rm lk}_{N_{17}}(3) \cong {\rm
lk}_{N_{17}}(8) \cong R_4$, (iii) ${\rm lk}_{N_{18}}(3) \cong
{\rm lk}_{N_{18}}(8) \cong R_3$ and $G_6(N_{18}) = (V, \{73, 31,
18, 84\})$, (iv) ${\rm lk}_{N_{19}}(3) \cong {\rm lk}_{N_{19}}(8)
\cong R_3$ and $G_6(N_{19}) = (V, \{63, 31, 18, 86\})$, (v) ${\rm
lk}_{N_{20}}(3) \cong {\rm lk}_{N_{20}}(8) \cong R_3$ and
$G_6(N_{20})$ $= (V, \{74, 28, 83, 31\})$, (vi) ${\rm
lk}_{N_{21}}(3) \cong R_2$, ${\rm lk}_{N_{21}}(8) \cong R_3$ and
$G_6(N_{21}) = (V, \{48, 83, 37, 36\})$, (vii) ${\rm
lk}_{N_{22}}(3) \cong R_2$, ${\rm lk}_{N_{22}}(8) \cong R_3$ and
$G_6(N_{22}) = (V, \{28, 86, 63, 37, 38\})$, (viii) ${\rm
lk}_{N_{23}}(3) \cong R_1$ and ${\rm lk}_{N_{23}}(8) \cong R_3$,
(ix) ${\rm lk}_{N_{24}}(3) \cong {\rm lk}_{N_{24}}(8) \cong R_1$.
These imply that there is no isomorphism between any two members
of ${\cal N}_1$.

Observe that (i) ${\rm lk}_{N_{25}}(3) \cong R_3$ and ${\rm
lk}_{N_{25}}(8) \cong R_4$, (ii) ${\rm lk}_{N_{26}}(3) \cong {\rm
lk}_{N_{26}}(8) \cong R_3$ and $G_6(N_{26}) = (V, \{53, 38,
84\})$, (iii) ${\rm lk}_{N_{27}}(3) \cong {\rm lk}_{N_{27}}(8)
\cong R_3$, $G_6(N_{27}) = (V, \{78, 81, 13, 37\})$ and ${\rm
NEG}(N_{27})= \{24, 56\}$, (iv) ${\rm lk}_{N_{28}}(3) \cong {\rm
lk}_{N_{28}}(8) \cong R_3$, $G_6(N_{28})$ $= (V, \{18, 84, 43,
31\})$ and ${\rm NEG}(N_{28})= \{75, 56\}$, (v) ${\rm
lk}_{N_{29}}(3) \cong R_3$ and ${\rm lk}_{N_{29}}(8) \cong R_2$,
(vi)  ${\rm lk}_{N_{30}}(3) \cong R_1$ and ${\rm lk}_{N_{30}}(8)
\cong R_3$, (vii) ${\rm lk}_{N_{31}}(3) \cong {\rm lk}_{N_{31}}(8)
\cong R_2$. These imply that there is no isomorphism between any
two members of ${\cal N}_2$.

Observe that (i) ${\rm lk}_{N_{32}}(3) \cong {\rm lk}_{N_{32}}(8)
\cong R_3$, (ii) ${\rm lk}_{N_{33}}(3) \cong {\rm lk}_{N_{33}}(8)
\cong R_4$, (iii) ${\rm lk}_{N_{34}}(3) \cong {\rm
lk}_{N_{34}}(8) \cong R_2$, (iv)  ${\rm lk}_{N_{35}}(3) \cong
R_2$ and ${\rm lk}_{N_{35}}(8) \cong R_1$. These imply that there
is no isomorphism between any two members of ${\cal N}_3$.

Since a member of ${\cal N}_i$ is non-isomorphic to a member of
${\cal N}_j$ for $i \neq j$, the above imply Part $(a)$. Part
$(b)$ follows from the definition of $N_{k}$ for $8 \leq k\leq
35$. \hfill $\Box$

\bigskip

The 3-dimensional {\em Kummer variety} $K^3$ is the torus
$S^1\times S^1 \times S^1$ modulo the involution $\sigma : x
\mapsto -x$. It has 8 singular points corresponding to 8 elements
of order 2 in the abelian group $S^1\times S^1 \times S^1$. In
\cite{k1}, K\"{u}hnel showed that $N_5$ triangulates $K^3$. For a
topological space $X$, $C(X)$ denotes a cone with base $X$. Let
$H = D^2\times S^1$ denote the solid torus. As a consequence of
the above lemmas we get.

\begin{cor}$\!\!\!${\bf .} \label{topology}
All the $8$-vertex normal $3$-pseudomanifolds triangulate seven
distinct topological spaces, namely, $|S^3_{8, j}| = S^3$ for
$1\leq j \leq 38$, $|N_1|$, $|N_2| = S(S^1\times S^1)$, $|N_3|$,
$|N_4| = H \cup (C(\partial H))$, $|N_5| = |N_6| = K^3$ and $|N_i|
= S(\RR P^2)$ for $7 \leq i\leq 35$.
\end{cor}

\noindent {\bf Proof.} Let $K$ be a $8$-vertex normal
$3$-pseudomanifold. If $K$ is a combinatorial $3$-sphere then it
triangulates the 3-sphere $S^3$.

If $K$ is not a combinatorial $3$-sphere then, by Lemma
\ref{le3.5} $(b)$, $|K|$ is (pl) homeomorphic to $|N_1|, \dots,
|N_6|$ or $|N_7|$. Since $N_2 = \Sigma_{78}T$, $|N_2|$ is
homeomorphic to the suspension $S(S^1\times S^1)$. In $N_{4}$,
the facets not containing the vertex 8 form a solid torus whose
boundary is the link of 8. This implies that $|N_4| = H \cup
(C(\partial H))$. It follows from Lemma \ref{le3.3} $(c)$ that
$|N_{6}|$ is (pl) homeomorphic to $|N_{5}|=K^3$. Since $N_{24}$
is isomorphic to the suspension $S^0_2 \ast R_{1}$, $|N_{24}| =
S(\RR P^2)$. Therefore, by Lemma \ref{le3.5} $(b)$, $|N_i|$ is
(pl) homeomorphic to $|N_{24}| = S(\RR P^2)$ for $7\leq i\leq
35$. The result now follows from Lemma \ref{le3.3} $(a)$. \hfill
$\Box$

\bigskip


A 3-dimensional {\em pseudo-complex} $K$ is an ordered pair
$(\Delta, \Phi)$, where $\Delta$ is a finite collection of
disjoint tetrahedra and $\Phi$ is a family of affine isomorphisms
between pairs of 2-faces of the tetrahedra in $\Delta$. Let $|K|$
denote the quotient space obtained from the disjoint union
$\sqcup_{\sigma\in \Delta}\sigma$ by setting $x = \varphi(x)$ for
$\varphi\in \Phi$. The quotient of a tetrahedron $\sigma \in
\Delta$ in $|K|$ is called a {\em $3$-simplex} in $|K|$ and is
denoted by $|\sigma|$. Similarly, the quotient of 2-faces, edges
and vertices of tetrahedra are called {\em $2$-simplices}, {\em
edges} and {\em vertices} in $|K|$ respectively. If $|K|$ is
homeomorphic to a topological space $X$ then $K$ is called a {\em
pseudo-triangulation} of $X$.
A 3-dimensional pseudo-complex $K = (\Delta, \Phi)$ is said to be
{\em regular} if (i) each 3-simplex in $|K|$ has four distinct
vertices and (ii) for $2 \leq i \leq 3$, no two distinct
$i$-simplices in $|K|$ have the same set of vertices. So, for $2
\leq i \leq 3$, an $i$-simplex $\alpha$ in $|K|$ is uniquely
determined by its vertices and denoted by $u_1\cdots u_{i+1}$,
where $u_1, \dots, u_{i+1}$ are vertices of $\alpha$. (But, the
edges in $|K|$ may not form a simple graph.) So, we can identify
a regular pseudo-complex $K = (\Delta, \Phi)$ with ${\cal K} :=
\{|\sigma| : \sigma \in \Delta\}$. Simplices and edges in $|K|$
are said to be  simplices and edges of ${\cal K}$. Clearly, a
pure 3-dimensional simplicial complex is a regular pseudo-complex.

Let ${\cal M}$ be a regular pseudo-triangulation of $X$ and
$abcd$, $abce$ be two 3-simplices in ${\cal M}$. If $ade, bde,
cde$ are not 2-simplices in ${\cal M}$, then ${\cal N} := ({\cal
M} \setminus \{abcd, abce\}) \cup \{abde, acde, bcde\}$ is also a
regular pseudo-triangulation of $X$. We say that ${\cal N}$ is
obtained from ${\cal M}$ by the {\em generalized bistellar
$1$-move} $\kappa_{abc}$. If there is no edge between $d$ and $e$
in ${\cal M}$ then $\kappa_{F}$ is called a {\em bistellar
$1$-move}. If there exist 3-simplices of the form $xyuv$, $xzuv$,
$yzuv$ in a regular pseudo-triangulation ${\cal P}$ of $Y$ and
$xyz$ is not a 2-simplex, then ${\cal Q} := ({\cal P} \setminus
\{xyuv, xzuv, yzuv\}) \cup \{xyzu, xyzv\}$ is also a regular
pseudo-triangulation of $Y$. We say that ${\cal Q}$ is obtained
from ${\cal P}$ by the {\em generalized bistellar $2$-move}
$\kappa_{E}$, where $E$ is the common edge in $xyuv$, $xzuv$ and
$yzuv$. If $E$ is the only edge between $u$ and $v$ in ${\cal P}$
then $\kappa_{E}$ is called a {\em bistellar $2$-move}.

Let $M$ be a pseudo-triangulation of a closed 3-manifold and $N$
a 3-pseudomanifold. A simplicial map $f \colon M \to N$ is said
to be a {\em $k$-fold branched covering} (with discrete branch
locus) if there exists $U \subseteq V(N)$ such that
$|f||_{|M|\setminus f^{- 1}(U)} \colon |M| \setminus f^{- 1}(U)
\to |N| \setminus U$ is a $k$-fold covering. The smallest such
$U$ (so that $|f||_{|M|\setminus f^{- 1}(U)} \colon |M| \setminus
f^{- 1}(U) \to |N|\setminus U$ is a covering) is called the {\em
branch locus}. It is known that $N_1$ can be regarded as a
branched quotient of a regular hyperbolic tessellation (cf.
\cite{k3}). In \cite{k1}, K\"{u}hnel has shown that $N_5$ is a
2-fold branched quotient of a pseudo-triangulation of the
3-dimensional torus. Here we prove\,:

\begin{theo}$\!\!\!${\bf .} \label{bcover2}
$(a)$ $N_{24}$ is a $2$-fold branched quotient of a $14$-vertex
combinatorial $3$-sphere. \newline $(b)$ For $7 \leq i \leq 35$,
$N_i$ is a $2$-fold branched quotient of a $14$-vertex regular
pseudo-triangulation of the $3$-sphere.
\end{theo}


\begin{lemma}$\!\!\!${\bf .} \label{le3.6}
Let $M$ be a regular pseudo-triangulation of a $3$-manifold and
$N$ be a normal $3$-pseudomanifold. Let $f \colon M \to N$ be a
$k$-fold branched covering with at most two vertices in the
branch locus. If $\kappa_e \colon N \mapsto \widetilde{N}$ is a
bistellar $2$-move, then there exist $k$ generalized bistellar
$2$-moves $\kappa_{e_1}, \dots, \kappa_{e_k}$ such that
$\kappa_{e_k}(\cdots (\kappa_{e_1}(M)))$ is a $k$-fold branched
cover of $\widetilde{N}$.
\end{lemma}

\noindent {\bf Proof.} Let ${\rm lk}_N(e) = S^{1}_{3}(\{x, y,
z\})$. Let $f^{-1}(e)$ consists of the edges $e_1, \dots, e_k$.
Let the end points of $e_i$ be $u_i$, $v_i$, the 3-simplices
containing $e_i$ be $u_iv_ix_iy_i$, $u_iv_ix_iz_i$,
$u_iv_iy_iz_i$ and $f(x_i) = x$, $f(y_i) = y$, $f(z_i) = z$ for
$1 \leq i\leq k$. Since $xyz$ is not a simplex in $N$, it follows
that $x_iy_iz_i$ is not a 2-simplex in $M$. Let $M_i$ be the
pseudo-complex consists of $u_iv_ix_iy_i$, $u_i v_i x_i z_i$ and
$u_iv_iy_iz_i$. Since the number of vertices in the branched
locus is at most 2, it follows that the number of vertices common
in $M_i$ and $M_j$ is at most 2 for $i\neq j$. In particular,
$\#(\{x_i, y_i, z_i\} \cap \{x_j, y_j, z_j\}) \leq 2$. Therefore,
$x_jy_jz_j$ is not a 2-simplex in $\kappa_{e_i}(M)$. So, we can
perform generalized bistellar 2-move $\kappa_{e_j}$ on
$\kappa_{e_i}(M) = (M\setminus M_i) \cup \{x_iy_iz_iu_i, x_iy_iz_i
v_i\}$ for $i \neq j$. Clearly, $\widetilde{M} := \kappa_{e_k}(
\cdots \kappa_{e_1}(M))$ is a $k$-fold branched cover of
$\widetilde{N}$ (via the map $\tilde{f}$, where $\tilde{f}(w) =
f(w)$ for $w\in V(\widetilde{M}) = V(M)$ and $\tilde{f}(x_iy_iz_i
u_i) = xyzu$ and $\tilde{f}(x_iy_iz_iv_i) = xyzv$).  \hfill $\Box$

\begin{lemma}$\!\!\!${\bf .} \label{le3.7}
Let $M$ be a regular pseudo-triangulation of a $3$-manifold and
$N$ be a normal $3$-pseudomanifold. Let $f \colon M \to N$ be a
$k$-fold branched covering with at most two vertices in the
branch locus. If $\kappa_F \colon N \mapsto \widetilde{N}$ is a
bistellar $1$-move, then there exist $k$ generalized bistellar
$1$-moves $\kappa_{F_1}, \dots, \kappa_{F_k}$ such that
$\kappa_{F_k}(\cdots (\kappa_{F_1}(M)))$ is a $k$-fold branched
cover of $\widetilde{N}$.
\end{lemma}

\noindent {\bf Proof.} Let $F = xyz$ and ${\rm lk}_N(F)=\{u, v\}$.
Let $f^{-1}(F)$ consist of the 2-simplices $F_1, \dots, F_k$. Let
$F_i = x_iy_iz_i$ and the 3-simplices containing $F_i$ be
$x_iy_iz_iu_i$ and $x_iy_iz_iv_i$ and $f(x_i, y_i, z_i, u_i, $ $
v_i) = (x, y, z, u, v)$ for $1 \leq i\leq k$. Since $f$ is
simplicial, it follows that $x_iu_iv_i$, $y_iu_iv_i$ and
$z_iu_iv_i$ are not 2-simplices in $M$. Let $M_i$ be
pseudo-complex $\{x_iy_iz_iu_i, x_iy_iz_iv_i\}$. Since the number
of vertices in the branched locus is at most 2, it follows that
$x_ju_jv_j$, $y_ju_jv_j$ and $z_ju_jv_j$ are not 2-simplices in
$\kappa_{F_i}(M)$ for $i\neq j$. Then (by the similar arguments
as in the proof of Lemma \ref{le3.6}) $\kappa_{F_k}(\cdots
\kappa_{F_1}(M))$ is a $k$-fold branched cover of
$\widetilde{N}$. \hfill $\Box$

\bigskip

\noindent {\bf Proof of Theorem \ref{bcover2}.} If ${\cal I}$
denotes the boundary of the icosahedron then there exists a
simplicial 2-fold covering $f \colon {\cal I} \to R_1$. Consider
the simplicial map $\tilde{f} \colon S^0_2(\{a, b\}) \ast {\cal
I} \to S^0_2(\{c, d\}) \ast R_1$ given by $\tilde{f}(a) = c$,
$\tilde{f}(b) = d$ and $\tilde{f}(u) = f(u)$ for $u\in V({\cal
I})$. Then $\tilde{f}$ is a 2-fold branched covering with branch
locus $\{c, d\}$. Since $N_{24}$ is isomorphic to the suspension
$S^0_2 \ast R_{1}$, it follows that $N_{24}$ is a 2-fold branched
quotient of the 14-vertex combinatorial 3-sphere $S^0_2(\{a, b\})
\ast {\cal I}$ (with branch locus $\{3, 8\}$). This proves $(a)$.

The result now follows from Lemmas \ref{le3.5} $(a)$, \ref{le3.6}
and \ref{le3.7}. (In fact, to obtain a 2-fold branched cover
$\widetilde{N}_{14}$ of $N_{14}$ from $R_1\ast S^0_2$ one needs
one bistellar 1-move and then one generalized bistellar 1-move.
And all other moves required in the proof are bistellar moves on
regular pseudo-triangulations of $S^3$.) \hfill $\Box$


\begin{remark}$\!\!\!${\bf .} \label{remark}
{\rm The combinatorial $3$-sphere $R_1 \ast S^0_2$ is a $2$-fold
branched cover of $N_{24}$ and $N_{14}$ can be obtained from
$N_{24}$ by a bistellar 1-move. Now, if $f \colon M \to N_{14}$
is a $2$-fold branched covering and $M$ is a combinatorial
$3$-manifold then (since ${\rm lk}_{N_{14}}(8)$ is a $7$-vertex
triangulated $\RR P^2$) the link of any vertex in $f^{- 1}(8)$ is
a 14-vertex triangulated $S^2$ and hence $f_0(M) > 14$.
(Similarly, for $i \neq 24$, if $N_i$ is a branched quotient of a
combinatorial 3-manifold $M$ then $f_0(M) > 14$.) So, there does
not exist a combinatorial $3$-sphere $M$ which is a branched
cover of $N_{14}$ and which can be obtained from $R_1 \ast S^0_2$
by proper bistellar moves. }
\end{remark}

In \cite{a2}, Altshuler observed that $N_1$ is orientable and
$|N_1|$ is simply-connected. In \cite{l1}, Lutz showed that
$(H_1(N_1), H_2(N_1), H_3(N_1)) = (0, \ZZ^8, \ZZ)$. The normal
$3$-pseudomanifold $N_3$ is the only among all the 35 which has
singular vertices of different types, namely one singular vertex
whose link is a triangulated torus and four singular vertices
whose links are triangulated real projective planes. Using
polymake (\cite{gj}), we find that $(H_1(N_3), H_2(N_3),
H_3(N_3))$ $ = (0, \ZZ^2\!\oplus \ZZ_2, 0)$. We summarized all the
findings about $N_1, \dots, N_{35}$ in the next table.

\bigskip


\noindent {\bf Table 1\,:} 8-vertex normal 3-pseudomanifolds which
are not combinatorial 3-manifolds.


\begin{center}
{\small
\begin{tabular}{|c|c|c|c|c|l|} \hline
&&&&&  \\[-3.5mm]
$X$ & $f$-vector & $\chi(X)$& $\!n_s(X)\!$ & links of sing- & Geometric carriers,  \\
&$ (f_1, f_2, f_3)$&&& ular vertices & Homology $(H_1, H_2, H_3)$ \\[1mm]\hline
&&&&& \\[-3.5mm]
$N_1$& $(28, 56, 28)$& 8& 8& all are $T$  & $|N_1|$ is simply connected,  \\
&&&&& $(H_1, H_2, H_3) = (0, \ZZ^8, \ZZ)$ \\[1mm]\hline
&&&&& \\[-3.5mm]
$N_2$& $(28, 44, 22)$& 2& 2& both are $T$ & $|N_2| = S(S^1 \times
S^1)$
\\[1mm] \hline
&&&&& \\[-3.5mm]
$N_3$& $(28, 46, 23)$& 3& 5& $T, R_2, R_2, $ & $ (H_1, H_2, H_3) =$ \\
&& & & $R_3, R_3$ & $~~~~~~~ (0, \ZZ^2\!\oplus \ZZ_2, 0)$ \\[1mm]\hline
&&&&& \\[-3.5mm]
$N_4$& $(28, 42, 21)$& 1& 1& $T$ & $|N_4| = H \cup (C(\partial H))$\\[1mm]\hline
&&&&& \\[-3.5mm]
$N_5$& $(28, 48, 24)$& 4& 8& all are $R_4$ & $|N_5| = K^3$\\[1mm]\hline
&&&&& \\[-3.5mm]
$N_6$& ,,& ,,& ,,& all are $R_3$ & $|N_6| = K^3$\\[1mm]\hline
&&&&& \\[-3.5mm]
$N_7$& $(28, 42, 21)$& 1& 2& both are $R_4$ & $|N_7| = S(\RR P^2)$\\[1mm]\hline
&&&&& \\[-3.5mm]
$N_i$, $8\leq i \leq 15$& ,, & ,,& ,, & both are in  & $|N_i| = S(\RR P^2)$\\
& & & & $\{R_1, \dots, R_4\}$ & \\[1mm]\hline
&&&&& \\[-3.5mm]
$N_i$, $16\leq i \leq 24$& $(27, 40, 20 )$ & ,,& ,, & ,,  & ~~~~~~~~~~
,,\\[1mm]\hline
&&&&& \\[-3.5mm]
$N_i$, $25\leq i \leq 31$& $(26, 38, 19)$ & ,,& ,, & ,,  & ~~~~~~~~~~
,,\\[1mm]\hline
&&&&& \\[-3.5mm]
$N_i$, $32\leq i \leq 35$& $(25, 36, 18)$ & ,,& ,, & ,,  & ~~~~~~~~~~
,,\\[1mm]\hline
\multicolumn{6}{l}{} \\[-2mm]
\multicolumn{6}{l}{[Here $K^3$ is the 3-dimensional Kummer
variety, $H = D^2\times S^1$ is the solid torus, $S(Y)$ is }  \\
\multicolumn{6}{l}{the topological suspension of $Y$ and $n_s(X)$
is the number of singular vertices in $X$.]}
\end{tabular}
}
\end{center}

\newpage

\begin{eg}$\!\!\!${\rm {\bf :}} \label{M310}
{\rm For $d\geq 2$, let
$$
K^{d}_{2d+3} = \{v_{i}\cdots v_{j-1}v_{j+1}\cdots v_{i+d+1} :
i+1\leq j\leq i+d, 1\leq i\leq 2d+3\}
$$
(additions in the suffixes are modulo $2d+3$). It was shown in
\cite{k2} the following\,: (i) $K^{\hspace{.2mm}d}_{2d+3}$ is a
triangulated $d$-manifold for all $d\geq 2$, (ii)
$K^{\hspace{.2mm}d}_{2d + 3}$ triangulates $S^{\hspace{.2mm}d -
1} \times S^{\hspace{.2mm}1}$ for $d$ even, and triangulates the
twisted product $\TPSSD$ (the twisted $S^{\hspace{.2mm}d -
1}$-bundle over $S^{\hspace{.2mm}1}$) for $d$ odd. For $d\geq 3$,
$K^{d}_{2d+3}$ is the unique non-simply connected $(2d+3)$-vertex
triangulated $d$-manifold (cf. \cite{bd8}). The combinatorial
3-manifolds $K^{\hspace{.2mm}3}_9$ was first constructed by
Walkup in \cite{w}.

From $K^3_9$ we construct the following 10-vertex combinatorial
3-manifold.
\begin{eqnarray*}
A^{3}_{10} & := & (K^3_9 \setminus\{v_1v_2v_3v_5, v_2v_3v_5v_6,
v_3v_5v_6v_7, v_3v_4v_6v_7, v_4v_6v_7v_8\}) \cup \\
&& \quad \{v_0v_1v_2v_3, v_0v_1v_2v_5, v_0v_1v_3v_5, v_0v_2v_3v_6,
v_0v_2v_5v_6, v_0v_3v_5v_7, \\
&& \qquad  v_0v_5v_6v_7, v_0v_3v_4v_6, v_0v_3v_4v_7, v_0v_4v_6v_8,
v_0v_4v_7v_8, v_0v_6v_7v_8\}.
\end{eqnarray*}
[Geometrically, first we remove a pl 3-ball consisting of five
3-simplices from $|K^3_9|$. This gives a pl 3-manifold with
boundary and the boundary is a 2-sphere. Then we add a cone with
base this boundary and vertex $v_0$. So, the new polyhedron
$|A^{3}_{10}|$ is pl homeomorphic to $|K^3_9|$. This implies that
the simplicial complex $A^{3}_{10}$ is a combinatorial
3-manifold.]

The only non-edge in $A^{3}_{10}$ is $v_0v_9$ and there is no
common $2$-face in the links of $v_0$ and $v_9$ in $A^{3}_{10}$.
So, $A^{3}_{10}$ does not allow any bistellar $1$-move. So,
$A^{3}_{10}$ is a 10-vertex non-neighbourly combinatorial
3-manifold which does not admit any bistellar 1-move.

Similarly, from $K^4_{11}$ we construct the following 12-vertex
triangulated 4-manifold.
\begin{eqnarray*}
A^{4}_{12} \!\!& := &\!\! (K^4_{11} \setminus\{v_1v_2v_3v_4v_6,
v_2v_3v_4v_6v_7,
v_3v_4v_6v_7v_8, v_4v_6v_7v_8v_9, v_4v_5v_7v_8v_9, v_5v_7v_8v_9v_{10}\})  \\
&& \quad \cup \{v_0v_1v_2v_3v_4, v_0v_1v_2v_3v_6, v_0v_1v_2v_4v_6,
v_0v_1v_3v_4v_6, v_0v_2v_3v_4v_7, v_0v_2v_3v_6v_7, \\
&& \qquad v_0v_2v_4v_6v_7, v_0v_3v_4v_6v_8, v_0v_3v_4v_7v_8,
v_0v_3v_6v_7v_8, v_0v_4v_6v_7v_9, v_0v_4v_6v_8v_9, \\
&& \qquad v_0v_4v_7v_8v_9, v_0v_4v_5v_7v_9, v_0v_4v_5v_8v_9,
v_0v_4v_7v_8v_9, v_0v_5v_7v_8v_{10}, v_0v_5v_7v_9v_{10}, \\
&& \qquad v_0v_5v_8v_9v_{10}\}.
\end{eqnarray*}
The only non-edge in $A^{4}_{12}$ is $v_0v_{11}$ and there is no
common $2$-face in the links of $v_0$ and $v_{11}$ in
$A^{4}_{12}$. So, $A^{4}_{12}$ does not allow any bistellar
$1$-move. So, $A^{4}_{12}$ is a 12-vertex non-neighbourly
triangulated 4-manifold which does not admit any bistellar 1-move.

By the same way, one can construct a $(2d+4)$-vertex
non-neighbourly triangulated $d$-manifold $A^d_{2d+4}$ (from
$K^{d}_{2d+3}$) which does not admit any bistellar 1-move for all
$d\geq 3$. }
\end{eg}

\begin{eg}$\!\!\!${\rm {\bf :}} \label{N39}
{\rm Let $N_{3}$ be as in Example \ref{pm38}. Let $M$ be obtained
from $N_{3}$ by starring two vertices $u$ and $v$ in the facets
$1248$ and $3568$ respectively, i.e., $M = \kappa_{1248}
(\kappa_{3568}(N_{3}))$. Then $M$ is a 10-vertex normal
3-pseudomanifold. Let $B^3_9$ be obtained from $M$ by identifying
the vertices $u$ and $v$. Let the new vertex be 9. Then
\begin{eqnarray*}
B^3_9 & := & (N_{3} \setminus \{1248, 3568\}) \cup \{1249, 1289,
1489, 2489, 3569, 3589, 3689, 5689\}.
\end{eqnarray*}
The degree 3 edges in $B^3_9$ are $16$, $17$ and $67$. But, none
of these edges is removable. So, no bistellar $2$-moves are
possible from $B^3_9$. The only non-edge in $B^3_9$ is $79$.
Since there is no common $2$-face in the links of $7$ and $9$, no
bistellar 1-move is possible. So, $B^3_9$ is a $9$-vertex
non-neighbourly 3-pseudomanifold  which does not admit any proper
bistellar move. }
\end{eg}


\section{Proofs}

For $n\geq 4$, by an $S^{\hspace{.2mm}2}_n$ we mean a
combinatorial 2-sphere on $n$ vertices. If $\kappa_{\beta} \colon
M \mapsto N$ is a bistellar 1-move then $\deg_N(v) \geq
\deg_M(v)$ for $v\in V(M)$. Here we prove the following.

\begin{lemma}$\!\!\!${\bf .} \label{le4.1}
Let $M$ be an $n$-vertex $3$-pseudomanifold and $u$ be a vertex
of degree $4$. If $n\geq 6$ then there exists a bistellar
$1$-move $\kappa_{\beta} : M \mapsto N$ such that $\deg_N(u) = 5$.
\end{lemma}

\noindent {\bf Proof.} Let ${\rm lk}_M(u) = S^2_4(\{a, b, c,
d\})$ and $\beta = abc$. Let ${\rm lk}_M(\beta) = \{u, x\}$. If
$x = d$ then the induced complex $K = M[\{u, a, b, c, d\}]$ is a
3-pseudomanifold. Since $n \geq 6$, $K$ is a proper subcomplex of
$M$. This is not possible. So, $x \neq d$ and hence $ux$ is a
non-edge in $M$. Then $\kappa_{\beta}$ is a bistellar 1-move.
Since $ux$ is an edge in $\kappa_{\beta}(M)$, $\kappa_{\beta}$ is
a required bistellar 1-move. \hfill $\Box$

\begin{lemma}$\!\!\!${\bf .} \label{le4.2}
Let $M$ be an $n$-vertex $3$-pseudomanifold and $u$ be a vertex
of degree $5$. If $n \geq 7$ then there exists a bistellar
$1$-move $\kappa_{\beta} : M \mapsto N$ such that $\deg_N(u) = 6$.
\end{lemma}

\noindent {\bf Proof.} Since $\deg_M(u) = 5$, the link of $u$ in
$M$ is of the form $S^0_2(\{a, b\})\ast S^1_3(\{x, y, z\})$ for
some vertices $a, b, x, y, z$ of $M$. If both $xyza$ and $xuzb$
are facets then the induced subcomplex $M[\{x, y, z, u, a, b\}]$
is a 3-pseudomanifold. This is not possible since $n \geq 7$. So,
without loss of generality, assume that $xyza$ is not a facet.
Again, if $xyab, xzab$ and $yzab$ all are facets then the induced
subcomplex $M[\{u, x, y, z, a, b\}]$ is a 3-pseudomanifold, which
is not possible. So, assume that $xyab$ is not a facet.

Consider the face $\beta = xya$. Suppose ${\rm lk}_M(\beta) =
\{u, w\}$. From the above, $w\not\in \{z, b\}$. So, $uw$ is a
non-edge and hence $\kappa_{\beta}$ is a required bistellar
1-move. \hfill $\Box$

\begin{lemma}$\!\!\!${\bf .} \label{le4.3}
Let $M$ be a non-neighbourly $8$-vertex $3$-pseudomanifold and $u$
be a vertex of degree $6$. If the degree of each vertex is at
least $6$, then there exists a bistellar $1$-move $\kappa_{\tau}
: M \mapsto N$ such that $\deg_N(u) = 7$.
\end{lemma}

\noindent {\bf Proof.} Let $u$ be a vertex with $\deg_M(u) =6$
and $uv$ be a non-edge. Let $L = {\rm lk}_M(u)$.

\smallskip

\noindent {\em Claim.} There exists a $2$-face $\tau$ such that
$\tau \cup \{u\}$ and $\tau \cup \{v\}$ are facets.

First consider the case when there exists a vertex $w$ such that
$\deg_L(w) = 5$. Let ${\rm lk}_L(w) (= {\rm lk}_M(uw)) = C_5(1,
2, 3, 4, 5)$.

Let $K = {\rm lk}_M(w)$. Since $\deg(v) = 6$, $vw$ is an edge.
Thus $K$ contains 7 vertices. If one of $12v, \dots, 45v, 51v$ is
a $2$-face, say $12v$, then $12wv$ and $12wu$ are facets. In this
case, $\tau = 12w$ serves the purpose.  So, assume that $12v,
\dots, 45v, 51v$ are non-faces in $K$. Then there are at least
three $2$-faces (not containing $u$) containing the edges $12,
\dots, 45, 51$ in $K$. Also, there are at least three $2$-faces
containing $v$ in $K$. So, the number of $2$-faces in $K$ is at
least $11$. This implies that $\deg_K(v) = 3$ or 4 and $K$ is a
7-vertex $\RR P^2$ or $P_4$. Since $\deg_K(u) = 5$, it follows
that $K$ is isomorphic to $R_2$, $R_3$ or $P_4$ (defined in
Section 2). In each case, (since $\deg_K(u) = 5$, $\deg_K(v) = 3$
or 4 and $uv$ is a non-edge) there exists an edge $\alpha$ in $K$
such that $\alpha \cup \{u\}$ and $\alpha \cup \{v\}$ are
$2$-faces in $K$ and hence $\tau = \alpha \cup \{w\}$ serves the
purpose.

Now, assume that $L$ has no vertex of degree 5. Then $L$ must be
of the form $S^0_2(\{a_1, a_2\}) \ast S^0_2(\{b_1, b_2\}) \ast
S^0_2(\{c_1, c_2\})$. If possible let $a_i b_j c_k v$ is not a
facet for $1 \leq i, j, k \leq 2$. Consider the $2$-face $a_1 b_1
c_1$. There exists a vertex $x \neq u$ such that $a_1 b_1 c_1 x$
is a facet. Assume, without loss of generality, that $a_1 b_1 c_1
a_2$ is a facet. Since $\deg(c_1) > 5$ (respectively, $\deg(b_1)
> 5$), $a_1 a_2 b_2 c_1$ (respectively, $a_1 a_2 b_1 c_2$) is not
a facet. So, the facet (other than $a_1 b_2 c_1u$) containing $a_1
b_2 c_1$ must be $a_1 b_2 c_1 c_2$. Similarly, the facet (other
than $a_1 b_1 c_2 u$) containing $a_1 b_1 c_2$ must be $a_1 b_1
b_2 c_2$. Then $a_1 b_2 c_1 c_2$, $a_1 b_1 b_2 c_2$ and $a_1 b_2
c_2 u$ are three facets containing $a_1 b_2 c_2$, a
contradiction. This proves the claim.

By the claim, there exists a 2-simplex $\tau$ such that ${\rm
lk}_M(\tau) = \{u, v\}$. Since $uv$ is a non-edge of $M$,
$\kappa_{\tau} : M \mapsto \kappa_{\tau}(M) = N$ is a bistellar
1-move. Since $uv$ is an edge in $N$, it follows that $\deg_N(u) =
7$. \hfill $\Box$

\bigskip

\noindent {\bf Proof of Theorem 1.} Let $M$ be an $8$-vertex
$3$-pseudomanifold. Then, by Lemma \ref{le4.1}, there exist
bistellar 1-moves $\kappa_{A_1}, \dots, \kappa_{A_k}$, for some
$k\geq 0$, such that the degree of each vertex in
$\kappa_{A_k}(\cdots(\kappa_{A_1}(M))$ is at least 5. Therefore,
by Lemma \ref{le4.2}, there exist bistellar 1-moves
$\kappa_{A_{k+1}}, \dots, \kappa_{A_l}$, for some $l \geq k$,
such that the degree of each vertex in $\kappa_{A_l}(\cdots
\kappa_{A_k}(\cdots (\kappa_{A_1}(M)))$ is at least 6. Then, by
Lemma \ref{le4.3}, there exist bistellar 1-moves
$\kappa_{A_{l+1}}, \dots, \kappa_{A_m}$, for some $m \geq l$,
such that the degree of each vertex in $\kappa_{A_m}(\cdots
\kappa_{A_l}(\cdots \kappa_{A_k}(\cdots (\kappa_{A_1}(M))))$ is
7. This proves the theorem. \hfill $\Box$

\begin{lemma}$\!\!\!${\bf .} \label{link=sphere}
Let $K$ be an $8$-vertex combinatorial $3$-manifold. If $K$ is
neighbourly then $K$ is isomorphic to $S^{3}_{8, 35}$, $S^{3}_{8,
36}$, $S^{3}_{8, 37}$ or $S^{3}_{8, 38}$.
\end{lemma}

\noindent {\bf Proof.} Since $K$ is a neighbourly combinatorial
3-manifold, by Proposition 2.1, the link of any vertex is
isomorphic to $S_5, \dots, S_8$ or $S_9$.

\smallskip

\noindent {\sf Claim}\,: The links of all the vertices can not be
isomorphic to $S_9$ ($= S^0_2\ast C_5$).

\smallskip

Otherwise, let ${\rm lk}(8) = S^0_2(6, 7) \ast C_5(1, 2, \dots,
5)$. Consider the vertex 2. Since the degree of 2 is 7,  $1267$ or
$2367$ is not a facet. Assume without loss of generality that
$1267$ is not a facet. Again, if $1236$ is a facet then
$\deg_{{\rm lk}(2)}(6) = 3$ and hence ${\rm lk}(2)\not\cong S_9$.
So, $1236$ is not a facet. Similarly, $1256$ is not a facet. Then
the facet other than $1268$ containing $126$ must be $1246$.
Similarly, $1247$ is a facet. This implies that ${\rm lk}(2) =
S^0_2(6, 7)\ast C_5(1, 4, 5, 3, 8)$. Thus $\deg(26) = 5$.
Similarly, $\deg(16) = \deg(36) = \deg(46) = \deg(56) = 5$. Then,
the 7-vertex 2-sphere ${\rm lk}(6)$ contains five vertices of
degree 5. This is not possible. This proves the claim.

\smallskip

\noindent {\bf Case 1.} Consider the case when $K$ has a vertex,
(say $8$) whose link is isomorphic to $S_8$. Assume, without loss
of generality, that the facets containing the vertex $8$ are
$1238$, $1268$, $1348$, $1458$, $1568$, $2348$, $2478$, $2678$,
$4578$ and $5678$. Since $\deg(3) =7$, $1234 \not\in K$. Hence the
facet other than $1238$ containing the face $123$ is one of
$1235$, $1236$ or $1237$.

If $1236 \in K$ then, clearly, $\deg(17) =3$ or $4$. If $\deg(17)
=4$, then on completing ${\rm lk}(1)$, we see that $1457$, $1567
\in K$, thereby showing that $\deg(5)=5$, an impossibility.
Hence, $\deg(17) =3$ and therefore $1457 \in K$. There are two
possibilities for the completion of ${\rm lk}(1)$. If $1347$,
$1356$, $1357 \in K$, from the links of $4$ and $3$, we see that
$2346$, $2467$, $3467$, $3567 \in K$. Here, $\deg(5) =6$. If
$1346$, $1467$, $1567 \in K$, then $\deg(5) =5$. Thus,
$1236\not\in K$.


\noindent {\bf Subcase 1.1.} $1235 \in K$. Since $\deg(1)=7$,
either $1345$ or $1256$ is a facet. In the first case, $1257$,
$1267$, $1567 \in K$. Here, $\deg(6)=5$, a contradiction. So,
$1256\in M$ and hence $1347$, $1357$, $1457 \in K$. From the links
of the vertices $1$, $4$, $7$ and $5$ we see that $1256$, $2346$,
$2467$, $3467$, $3567$, $2356 \in K$. Here, $K \cong S^{3}_{8,38}$
by the map $(1,5,8,6)(2,7)(3,4)$.


\noindent {\bf Subcase 1.2.} $1237 \in K$. By the same argument as
in Subcase 1.1 (replace the vertex 1 by vertex 2) we get $1267$,
$2345$, $2357$, $2457 \in K$. From ${\rm lk}(1)$ and ${\rm
lk}(7)$, $1346$, $1456$, $3456$, $1367$, $3567 \in K$. Here, $K
\cong S^{3}_{8, 38}$ by the map $(1,7,8,6)(2,5)(3,4)$.

\smallskip

\noindent {\bf Case 2.} $K$ has no vertex whose link is isomorphic
to $S_8$ but has a vertex whose link is isomorphic to $S_6$.
Using the same method as in Subcase 1.1, we find that $K \cong
S^{3}_{8, 37}$.

\smallskip

\noindent {\bf Case 3.} $K$ has no vertex whose link is isomorphic
to $S_8$ or $S_6$ but has a vertex whose link is isomorphic to
$S_7$. Using the same method as in Subcase 1.1, we find that $K
\cong S^{3}_{8, 36}$.

\smallskip

\noindent {\bf Case 4.} $K$ has no vertex whose link is isomorphic
to $S_6$, $S_7$ or $S_8$ but has a vertex (say 8) whose link is
isomorphic to $S_5$. The facets through $8$ can be assumed to be
$1238$, $1278$, $1348$, $1458$, $1568$, $1678$, $2348$, $2458$,
$2568$ and $2678$. Clearly, $1234$, $1267 \not\in K$. If
$\deg(15) = 6$, then from ${\rm lk}(1)$ and ${\rm lk}(5)$, we see
that $1235$, $1345$, $2345 \in K$, thereby showing that $\deg(3)
= 5$. Hence $1237 \in K$. Now, we can assume without loss of
generality that the facets required to complete ${\rm lk}(1)$ are
$1347$, $1457$ and $1567$. Now, consider ${\rm lk}(2)$. If
$\deg(27) = 6$, then after completing the links of 2 and 7 we
observe that $\deg(4) = 6$. Hence $\deg(23) = 6$. The links of
$2$, $7$ and $6$ show that $2345$, $2356$, $2367$, $3467$, $4567$
and $3456 \in K$. Here, $K \cong S^3_{8, 35}$ by the map $(2, 3,
4, 5, 6, 7, 8)$. This completes the proof. \hfill $\Box$

\begin{lemma}$\!\!\!${\bf .} \label{link=T}
Let $K$ be an $8$-vertex neighbourly normal $3$-pseudomanifold.
If $K$ has one vertex whose link is the $7$-vertex torus $T$ then
$K$ is isomorphic to $N_{1}$, $N_{2}$, $N_{3}$ or $N_{4}$.
\end{lemma}

\noindent {\bf Proof.} Let us assume that $V(K) = \{1, \dots, 8\}$
and the link of the vertex 8 is the 7-vertex torus $T$. So, the
facets containing 8 are $1248$, $1268$, $1348$, $1378$, $1568$,
$1578$, $2358$, $2378$, $2458$, $2678$, $3468$, $3568$, $4578$
and $4678$. We have the following cases.

\smallskip

\noindent {\bf Case 1.} There is a vertex (other than the vertex
8), say 7, whose link is isomorphic to $T$. Then ${\rm lk}(7)$
has no vertex of degree 3 and hence $2367, 1457, 1237, 1357
\not\in K$. This implies that the facet (other than $1378$)
containing $137$ is $1367$ or $1347$. In the first case, ${\rm
lk}(17) = C_6(5, 8, 3, 6, 4, 2)$. Thus, $1367, 1467, 1247, 1257
\in K$. Then, from the links of $67$ and $37$, we get $2567,
3567, 2347, 3457$ $\in K$. Now, from ${\rm lk}(34)$, $1346
\not\in K$. Then, from the links of $36, 34, 23, 14$ and $26$, we
get $1236, 2346, 1345, 1235, 1456, 2456 \in K$. Here, $K = N_{1}$.

In the second case, ${\rm lk}(37) = C_6(2, 8, 1, 4, 6, 5)$. Thus,
$1347, 3467, 3567, 2357 \in K$. Now, from the links of $47$ and
$67$, we get $1247, 2457, 1567, 1267 \in K$. Here, $K = N_{2}$.

\smallskip

\noindent {\bf Case 2.} There is a vertex whose link is a 7-vertex
$\RR P^2$.

\smallskip

\noindent {\sf Claim} \,: There exists a vertex in $K$ whose link
is isomorphic to $R_{2}$.

\smallskip

If there is vertex whose link is isomorphic to $R_2$ then we are
done. Otherwise, since ${\rm Aut}({\rm lk}(8))$ acts transitively
on $\{1, \dots, 7\}$, assume that ${\rm lk}(4) \cong R_{3}$
(respectively, $R_4$). Since $(1, 2, 5, 7, 6, 3) \in {\rm
Aut}({\rm lk}(8))$, we may assume that the degree 4 vertex
(respectively vertices) in ${\rm lk}(4)$ is 1 (respectively are
$1$, $5$, $6$). Then, from ${\rm lk}(4)$, $1247$, $1347$, $2467
\in K$. This implies that ${\rm lk}(7)$ is a non-sphere and
$\deg(67)=3$. Hence ${\rm lk}(7) \cong R_{2}$. This proves the
claim.

By the claim, we can assume that ${\rm lk}(4) \cong R_{2}$.
Again, we may assume that the vertex 1 is of degree 3 in ${\rm
lk}(4)$. Then, from ${\rm lk}(4)$, $1234$, $2347$, $2456$,
$2467$, $3456$, $3457 \in K$. Considering the links of the edges
$36$, $26$, $27$, $25$ and $13$, we get $1256$, $1235$, $1357 \in
K$. Here, $K = N_{3}$.

\smallskip

\noindent {\bf Case 3.} Only singular vertex in $K$ is $8$. So,
the link of each vertex (other than vertex 8) is an $S^2_7$ (a
7-vertex 2-sphere). Since $8$ is a degree 6 vertex in ${\rm
lk}(u)$, it follows that ${\rm lk}(u)$ is isomorphic to one of
$S_{5}$, $S_{6}$ or $S_{7}$ (defined in Example \ref{wpm27}) for
any vertex $u\neq 8$. If ${\rm lk}(1) \cong S_{5}$, then (since
$(3, 4, 2, 6, 5, 7) \in {\rm Aut}({\rm lk}(8))$), we may assume
that the other degree 6 vertex in ${\rm lk}(1)$ is 3. Then, from
the links of $1$ and $3$, $1348$, $1234$, $1346$ are facets
containing $134$, a contradiction. If ${\rm lk}(1) \cong S_{6}$
then (since ${\rm lk}(18) = C_6(3, 4, 2, 6, 5, 7)$), we may
assume that the degree 5 vertices in ${\rm lk}(1)$ are 2, 3, and
5. Then ${\rm lk}(3)$ can not be an $S^2_7$, a contradiction. So,
${\rm lk}(1) \cong S_{7}$. Since ${\rm Aut}({\rm lk}(8))$ acts
transitively on $\{1, \dots, 7\}$, it follows that the link of
each vertex is isomorphic to $S_7$.

Since ${\rm lk}(18) = C_6(3, 4, 2, 6, 5, 7)$ and $(3, 4, 2, 6, 5,
7) \in {\rm Aut}({\rm lk}(8))$, we may assume that the degree 5
vertices in ${\rm lk}(1)$ are 4 and 5. Since ${\rm lk}(4) \cong
S_7$, it follows that $1456 \not\in K$. Then, from ${\rm lk}(1)$,
$1245, 1256, 1347, 1457 \in K$. Now, from the links of $4$ and
$5$, we get $3467$, $2356 \in K$. Then, from ${\rm lk}(2)$, $2367
\in K$. Here $K = N_{4}$. This completes the proof. \hfill $\Box$

\begin{lemma}$\!\!\!${\bf .} \label{linknot=T}
Let $K$ be an $8$-vertex neighbourly normal $3$-pseudomanifold.
If $K$ is not a combinatorial $3$-manifold and has no vertex whose
link is isomorphic to the $7$-vertex torus $T$ then $K$ is
isomorphic to $N_{5}, \dots, N_{14}$ or $N_{15}$.
\end{lemma}

\noindent {\bf Proof.} Let $n_s$ be the number of singular
vertices in $K$. Since $K$ is neighbourly, by Proposition
\ref{w2mfd}, the link of any vertex is either a 7-vertex $\RR
P^2$ or a 7-vertex $S^2$. So, the number of facets through a
singular (respectively non-singular) vertex is 12 (respectively
10). Let $f_3$ be the number of facets of $K$. Consider the set $S
= \{(v, \sigma) ~ : ~ \sigma$ is a facet of $K$ and $v\in \sigma$
is a vertex$\}$. Then $f_3 \times 4 = \#(S) = n_s\times 12 +
(8-n_s) \times 10 = 80 +2n_s$. This implies $n_s$ is even. Since
$K$ is not a combinatorial 3-manifold, it follows that $n_s \neq
0$ and hence $n_s \geq 2$. So, $K$ has at least two vertices whose
links are isomorphic to $R_2$, $R_3$ or $R_4$.

\smallskip

\noindent {\bf Case 1.} There exist (at least) two vertices whose
links are isomorphic to $R_4$. Assume that ${\rm lk}_M(8) = R_4$.
Then $1258, 1268, 1358, 1378, 1468, 1478, 2368, 2378, 2458, 2478,
3458, 3468 \in K$. Since $(1, 3, 4)(5, 6, 7)$, $(1, 2)(3, 4)\in
{\rm Aut}({\rm lk}(8))$, we may assume that ${\rm lk}(3)$ or
${\rm lk}(7) \cong R_4$.

\smallskip

\noindent {\bf Subcase 1.1.} ${\rm lk}(7) \cong R_4$. Since ${\rm
lk}_{{\rm lk}(7)}(8) = C_4(1, 3, 2, 4)$, it follows that 1, 2, 3,
4 are degree 5 vertices in ${\rm lk}(7)$. Since $(3, 4)(5, 6) \in
{\rm Aut}({\rm lk}(8))$, assume without loss that $136, 145 \in
{\rm lk}(7)$. Then, from ${\rm lk}(7)$, we get $1257, 1267, 1367,
1457, 2357, 2467, 3457, 3467\in K$. This shows that ${\rm lk}(2)$
is an $\RR P^{2}_{7}$. Since $3457, 3458 \in K$, it follows that
$2345 \not\in K$. Then, from ${\rm lk}(2)$, $2356, 2456 \in K$.
Then, from the links of $3$ and $4$, $1356, 1456\in K$. Here $K =
N_{5}$.

\smallskip

\noindent {\bf Subcase 1.2.} ${\rm lk}(7) \not\cong R_4$. So,
${\rm lk}(3)\cong R_4$. Since ${\rm lk}_{{\rm lk}(3)}(8) = C_6(1,
7, 2, 6, 4, 5)$, the degree 4 vertices in ${\rm lk}(3)$ are
either $5, 6, 7$ or $1, 2, 4$. In the first case, on completion of
${\rm lk}(3)$, we observe that $56$, $67$, $57$ remain non-edges
in $K$. So, the degree 4 vertices in ${\rm lk}(3)$ are $1, 2$ and
$3$. Then $1356$, $1367$, $2356$, $2357$, $3457$ and $3467$ are
facets. Since ${\rm lk}(7) \not\cong R_4$ and $\deg(78) = 4$,
either ${\rm lk}(7) \cong R_3$ or ${\rm lk}(7)$ is an $S^2_7$. In
the former case, $2567$ is a facet. This is not possible from
${\rm lk}(25)$. So, ${\rm lk}(7)$ is an $S^2_7$. Then, from ${\rm
lk}(7)$, $1467, 2457\in K$. Now, from ${\rm lk}(1)$, $1256 \in
K$. Here, $K = N_{7}$.

\smallskip

\noindent {\bf Case 2.} Exactly one vertex whose link is
isomorphic to $R_4$ and there exists a vertex whose link is
isomorphic to $R_3$. Using the same method as in Case 1, we find
that $K \cong N_{8}$.

\smallskip

\noindent {\bf Case 3.} Exactly one vertex whose link is
isomorphic to $R_4$, there is no vertex whose link is isomorphic
to $R_3$ and there exists (at least) a vertex whose link is
isomorphic to $R_2$. Using the same method as in Case 1, we find
that $K \cong N_{9}$.

\smallskip

\noindent {\bf Case 4.} There is no vertex whose link is
isomorphic to $R_4$ and there exist (at least) two vertices whose
links are isomorphic to $R_3$. Assume that ${\rm lk}_K(8) = R_4$,
so that $\deg(78) = 4$. Using the same method as in Case 1, we
get the following: (i) if ${\rm lk}_K(7) \cong R_3$ then $K =
N_{6}$ and (ii) if ${\rm lk}_K(7) \not \cong R_3$ then $K$ is
isomorphic to $N_{10}$ or $N_{11}$.

\smallskip

\noindent {\bf Case 5.} There is no vertex whose link is
isomorphic to $R_4$, there exists exactly one vertex whose link
is isomorphic to $R_3$ and there exists (at least) a vertex whose
link is isomorphic to $R_2$. Using the same method as in Case 1,
we find that $K$ is isomorphic to $N_{12}$ or $N_{13}$.

\smallskip

\noindent {\bf Case 6.} There is no vertex whose link is
isomorphic to $R_4$ or $R_3$ and there exist (at least) two
vertices whose links are isomorphic to $R_2$. Using the same
method as in Case 1, we find that $K$ is isomorphic to $N_{14}$
or $N_{15}$. This completes the proof. \hfill $\Box$

\bigskip

\noindent {\bf Proof of Theorem 2.} Since $S^3_{8,m}$'s are
combinatorial 3-manifolds and $N_{n}$'s are not combinatorial
3-manifolds, $S^3_{8, m} \not\cong N_{n}$ for $35 \leq m \leq 38$,
$1 \leq n \leq 15$. Part $(a)$ now follows from Lemmas
\ref{le3.1}, \ref{le3.4}. Part $(b)$ follows from Lemmas
\ref{link=sphere}, \ref{link=T} and \ref{linknot=T}. \hfill $\Box$

\begin{lemma}$\!\!\!${\bf .} \label{S0toS1}
Let ${\cal S}_0, \dots, {\cal S}_6$ be as in the proof of Lemma
$\ref{le3.2}$. If a combinatorial $3$-manifold $K$ is obtained
from a member of ${\cal S}_j$ by a bistellar $2$-move then $K$ is
isomorphic to a member of ${\cal S}_{j+1}$ for $0\leq j\leq 5$.
Moreover, no bistellar $2$-move is possible from a member of
${\cal S}_6$.
\end{lemma}

\noindent {\bf Proof.} Recall that ${\cal S}_0 = \{S^3_{8,35},
S^3_{8,36}, S^3_{8,37}, S^3_{8,38}\}$. The removable edges in
$S^{3}_{8,37}$ are $13$, $16$, $17$, $24$, $27$, $35$, $46$, $48$
and $58$. Since $(1,4)(2,7)(3,8) \in {\rm Aut}(S^{3}_{8,37})$, up
to isomorphisms it is sufficient to consider the bistellar
$2$-moves $\kappa_{27}$, $\kappa_{24}$, $\kappa_{48}$,
$\kappa_{58}$ and $\kappa_{46}$ only. Here $S^3_{8, 33} :=
\kappa_{27}(S^3_{8, 37})$, $S^3_{8, 30} := \kappa_{24}(S^3_{8,
37})$, $S^3_{8, 32} := \kappa_{48}(S^3_{8, 37})$, $S^3_{8, 31} :=
\kappa_{58}(S^3_{8, 37})$ and $\kappa_{46}(S^3_{8, 37}) \cong
S^3_{8, 31}$ by the map $(1, 4, 5)(2, 7)(3, 6, 8)$.

The removable edges in $S^{3}_{8, 38}$ are $13$, $38$, $78$, $27$,
$25$, $15$ and $46$. Since $(1, 2, 8)(7, 3, 5)$, $(1, 2)(3, 7)(4,
6) \in {\rm Aut}(S^{3}_{8, 38})$, it is sufficient to consider the
bistellar 2-moves $\kappa_{46}$ and $\kappa_{78}$ only. Here
$S^{3}_{8, 39} := \kappa_{46}(S^{3}_{8, 36})$ and
$\kappa_{78}(S^{3}_{8, 38}) \cong S^{3}_{8, 32}$ by the map $(1,
7, 8, 4, 6)(2, 3)$.

The removable edges in $S^{3}_{8, 36}$ are $13, 35, 58, 68, 46,
24, 27$, $17$. Since $(1, 5, 6, 2)(3, 8, 4, 7)$ is an
automorphism of $S^{3}_{8, 36}$, it is sufficient to consider the
bistellar 2-moves $\kappa_{58}$ and $\kappa_{68}$ only. Here
$\kappa_{58}(S^{3}_{8, 36}) = S^{3}_{8, 31}$ and
$\kappa_{68}(S^{3}_{8, 36}) \cong S^{3}_{8, 30}$ by the map $(1,
6, 4, 8, 2, 5, 7, 3)$.

The removable edges in $S^{3}_{8, 35}$ are $13, 35, 57, 71, 24,
46, 68$ and $82$. Since $(1, 2, \dots, 8)$, $(1, 8)(2, 7)(3,
6)(4, 5) \in {\rm Aut}(S^{3}_{8, 35})$, it is sufficient to
consider the bistellar 2-moves $\kappa_{68}$ only. Here
$\kappa_{68}(S^{3}_{8, 35}) \cong S^{3}_{8, 30}$ by the map $(1,
7, 3)(2, 8, 4, 5, 6)$. This proves the result for $j=0$.

By the same arguments as in the case for $j=0$, one proves for
the cases for $1\leq j\leq 5$. We summarize these cases in the
diagram given below. Last part follows from the fact that none of
$S^3_{8, 1}$, $S^3_{8, 3}$ or $S^3_{8, 3}$ has any removable
edges. \hfill $\Box$

\setlength{\unitlength}{3mm}

\begin{picture}(45,36.7)(-2,3.5)


\thicklines

\put(12,38){\mbox{$S^3_{8,35}$}} \put(18,38){\mbox{$S^3_{8,36}$}}
\put(24,38){\mbox{$S^3_{8,37}$}} \put(30,38){\mbox{$S^3_{8,38}$}}

\put(12.2,37.5){\line(-1,-2){1.7}}
\put(18.5,37.5){\line(-1,-2){1.7}}
\put(18,37.5){\line(-2,-1){6.5}}
\put(24.5,37.5){\line(-1,-2){1.7}}
\put(24,37.5){\line(-2,-1){6.5}} \put(25.5,37.5){\line(1,-2){1.7}}
\put(23.5,37.5){\line(-3,-1){11.5}}
\put(30.5,37.5){\line(-1,-2){1.7}}
\put(31.8,37.5){\line(1,-2){1.7}}

\put(9,33){\mbox{$S^3_{8,30}$}} \put(15,33){\mbox{$S^3_{8,31}$}}
\put(21,33){\mbox{$S^3_{8,33}$}} \put(27,33){\mbox{$S^3_{8,32}$}}
\put(33,33){\mbox{$S^3_{8,39}$}}

\put(10.2,32){\line(4,-5){2.2}} \put(10.6,32){\line(5,-2){7}}
\put(11.3,32.05){\line(4,-1){12}}

\put(15,32.25){\line(-2,-1){6.5}} \put(16,32){\line(1,-1){2.5}}
\put(15.2,32){\line(-1,-2){1.3}} \put(16.5,32){\line(5,-2){7}}

\put(17.5,31.9){\line(4,-1){11.5}} \put(18,32){\line(6,-1){17}}
\put(19.1,30.1){\line(-4,1){8}}

\put(22.6,32.4){\line(2,-5){1.25}} \put(23,32.4){\line(2,-1){6.5}}
\put(23.3,32.4){\line(4,-1){11.5}}

\put(27,32.4){\line(-5,-1){18}} \put(28,32.4){\line(-4,-1){13.5}}
\put(34.5,32.2){\line(2,-3){1.7}}
\put(29.8,32.4){\line(2,-1){5.8}}

\put(6,28){\mbox{$S^3_{8,24}$}} \put(12,28){\mbox{$S^3_{8,23}$}}
\put(18,28){\mbox{$S^3_{8,28}$}} \put(24,28){\mbox{$S^3_{8,25}$}}
\put(30,28){\mbox{$S^3_{8,29}$}} \put(36,28){\mbox{$S^3_{8,27}$}}

\put(1.5,24.5){\line(3,2){4.5}} \put(2,24.5){\line(3,1){9.5}}
\put(2.5,24.2){\line(4,1){15}} \put(3,24){\line(5,1){20}}

\put(7.5,24.5){\line(3,2){4.5}} \put(8.5,24.2){\line(4,1){15}}

\put(14,24.5){\line(3,1){9.5}} \put(14.5,24){\line(4,1){15}}

\put(18,24){\line(-3,1){10}} \put(18.7,24.4){\line(-3,2){4.6}}
\put(19.5,24.5){\line(3,2){4.5}} \put(20.5,24){\line(4,1){15}}

\put(23.5,24){\line(-4,1){14.7}} \put(24.5,24.3){\line(0,1){3}}

\put(29.5,24){\line(-4,1){14.7}} \put(30,24.2){\line(-3,1){10.2}}
\put(30.5,24.2){\line(-3,2){4.8}}

\put(36,24.2){\line(-3,1){10}} \put(30.5,24.3){\line(0,1){3}}
\put(31.5,24.5){\line(3,2){4.5}}

\put(41.5,24){\line(-3,1){10}} \put(42.5,24.2){\line(-3,2){4.8}}

\put(0,23){\mbox{$S^3_{8,20}$}} \put(6,23){\mbox{$S^3_{8,11}$}}
\put(12,23){\mbox{$S^3_{8,18}$}} \put(18,23){\mbox{$S^3_{8,19}$}}
\put(24,23){\mbox{$S^3_{8,12}$}} \put(30,23){\mbox{$S^3_{8,22}$}}
\put(36,23){\mbox{$S^3_{8,21}$}} \put(42,23){\mbox{$S^3_{8,26}$}}

\put(0.5,19.3){\line(0,1){3}}

\put(6,19.2){\line(-3,2){4.5}} \put(6.5,19.3){\line(0,1){3}}
\put(7,19.2){\line(3,2){4.5}} \put(7.5,19.2){\line(3,1){10}}
\put(8.5,19.2){\line(4,1){15}} \put(8.5,18.8){\line(5,1){20}}

\put(12.5,19.3){\line(0,1){3}}

\put(17.8,19){\line(-5,1){15.8}} \put(18.5,19.2){\line(0,1){3.4}}
\put(18,19.2){\line(-3,2){4.5}} \put(19.5,19.2){\line(3,1){10}}

\put(23.2,18.5){\line(-5,1){20.5}}
\put(23.8,18.7){\line(-4,1){15.2}} \put(24.5,19.3){\line(0,1){3}}
\put(24,19.2){\line(-3,2){4.5}} \put(25.5,19.2){\line(3,1){10}}

\put(29.5,19){\line(-3,1){10}} \put(30,19.2){\line(-3,2){4.8}}

\put(35.5,19){\line(-6,1){21}} \put(36,19.3){\line(-5,1){16}}
\put(37,19.3){\line(0,1){3}} \put(36.3,19.4){\line(-3,2){4.5}}
\put(37.5,19.2){\line(3,2){4.5}}

\put(42.5,19.3){\line(0,1){3}}

\put(0,18){\mbox{$S^3_{8,14}$}} \put(6,18){\mbox{$S^3_{8,8}$}}
\put(12,18){\mbox{$S^3_{8,15}$}} \put(18,18){\mbox{$S^3_{8,16}$}}
\put(24,18){\mbox{$S^3_{8,9}$}} \put(30,18){\mbox{$S^3_{8,10}$}}
\put(36,18){\mbox{$S^3_{8,17}$}} \put(42,18){\mbox{$S^3_{8,34}$}}

\put(9,14){\line(-1,2){1.5}} \put(11,14){\line(3,1){12}}

\put(14,13.7){\line(-3,1){11.5}} \put(14.5,14){\line(-2,1){6}}
\put(15,14){\line(-1,2){1.65}} \put(15.5,14.2){\line(1,1){3}}
\put(17,14.2){\line(2,1){7}} \put(21,14){\line(-1,2){1.65}}

\put(26.5,13.7){\line(-5,1){18.2}} \put(27,14){\line(-2,1){7}}
\put(28,14.3){\line(1,1){3}}

\put(32.3,13.7){\line(-6,1){24}} \put(32.4,14){\line(-3,1){11.6}}
\put(33.2,14.2){\line(-2,1){6.9}} \put(34.2,14.2){\line(1,1){3}}

\put(9,13){\mbox{$S^3_{8,6}$}} \put(15,13){\mbox{$S^3_{8,5}$}}
\put(21,13){\mbox{$S^3_{8,13}$}} \put(27,13){\mbox{$S^3_{8,4}$}}
\put(33,13){\mbox{$S^3_{8,7}$}}

\put(14.8,9.2){\line(-3,2){4.5}} \put(15.5,9.3){\line(0,1){3}}
\put(16.4,9.2){\line(5,1){16.5}}

\put(20.8,9.2){\line(-3,2){4.5}} \put(22.3,9.3){\line(3,2){4.5}}
\put(21.5,9.3){\line(0,1){3}}

\put(27.5,9.3){\line(0,1){3}} \put(28.7,9.3){\line(3,2){4.5}}

\put(15,8){\mbox{$S^3_{8,2}$}} \put(21,8){\mbox{$S^3_{8,3}$}}
\put(27,8){\mbox{$S^3_{8,1}$}}

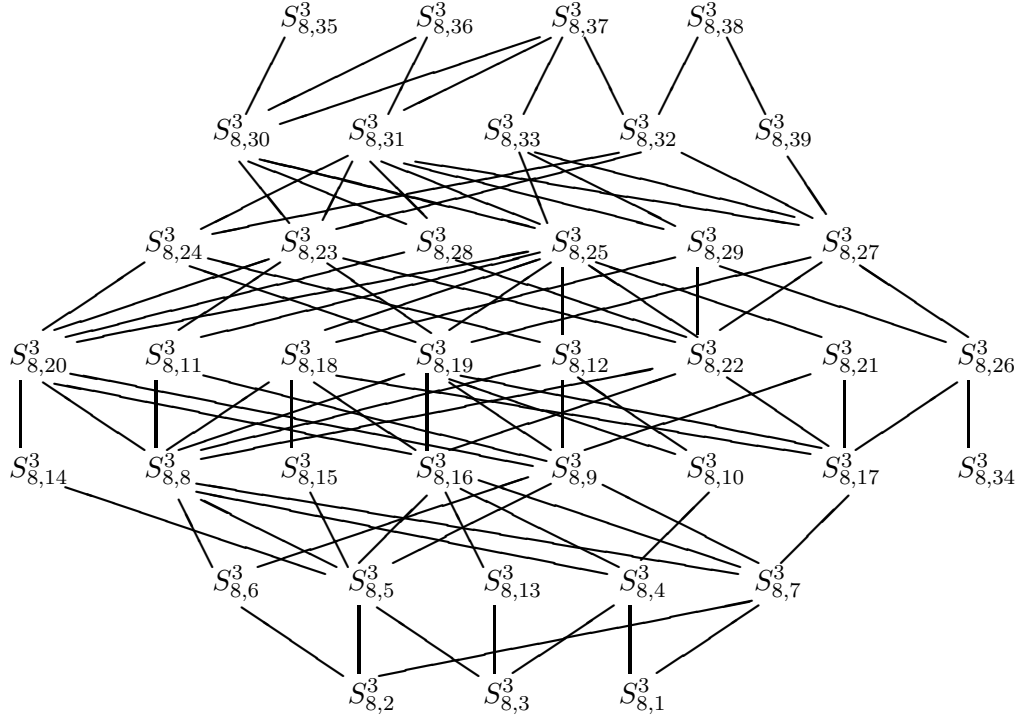
\put(-1,5.5){\mbox{{\bf Diagram 1\,:} Hasse diagram of the poset
of the 8-vertex combinatorial 3-manifolds}}

\put(8,3.8){\mbox{(the partial order relation is as defined in
Section 2)}}

\end{picture}

\begin{lemma}$\!\!\!${\bf .} \label{N0toN1}
Let ${\cal N}_0, \dots, {\cal N}_3$ be as in the proof of Lemma
$\ref{le3.5}$. If a $3$-pseudomanifold $K$ is obtained from a
member of ${\cal N}_{j}$ by a bistellar $2$-move then $K$ is
isomorphic to a member of ${\cal N}_{j+1}$ for $0\leq j\leq 2$.
Moreover, no bistellar $2$-move is possible from a member of
${\cal N}_3$.
\end{lemma}

\noindent {\bf Proof.} Recall that ${\cal N}_0=\{N_{1}, \dots,
N_{15}\}$. Since there are no degree 3 edges in $N_{1}$, $N_{2}$,
$N_{5}$ and $N_{6}$, no bistellar $2$-moves are possible from
$N_{1}$, $N_{5}$, $N_{6}$ or $N_{2}$. The degree 3 edges in
$N_{3}$ (respectively, in $N_{4}$) are $14, 16, 17, 36, 67$
(respectively, $13, 35, 57, 72, 24, 46, 61$). But, none of these
edges is removable. So, bistellar $2$-moves are not possible from
$N_{3}$ or $N_{4}$.

The removable edges in $N_{7}$ are $12, 14, 24, 56, 57$ and $67$.
Since $(1, 2)(6, 7)$, $(1, 2)(5, 6)$ and $(1, 5)(2, 6)(3, 8)(4,
7)$ are automorphisms of $N_{7}$, it follows that up to
isomorphisms, we only have to consider the bistellar 2-move
$\kappa_{67}$. Here, $N_{16} = \kappa_{67}(N_{7})$.

The removable edges in $N_{8}$ are $15, 17, 24, 56, 57$ and $67$.
Since $(1, 6)(2, 4)$, $(1, 6)(5, 7)$, $(2, 4)(5, 7) \in {\rm
Aut}(N_{8})$, we only consider the bistellar 2-moves
$\kappa_{24}$, $\kappa_{56}$ and $\kappa_{57}$. Here, $N_{17} =
\kappa_{24}(N_{8})$, $N_{18} = \kappa_{56}(N_{8})$ and $N_{19} =
\kappa_{57}(N_{8})$.

The removable edges in $N_{9}$ are $12, 23, 24$ and $67$. Since
$(1, 4)(6, 7) \in {\rm Aut}(N_{9})$, we consider only
$\kappa_{12}, \kappa_{23}$ and $\kappa_{67}$. Here, $N_{21} =
\kappa_{12}(N_{9})$, $N_{23} = \kappa_{23}(N_{9})$ and
$\kappa_{67}(N_{9}) = N_{16}$.

The removable edges in $N_{10}$ are $12, 14, 24, 56, 57$ and
$67$. Since $(1,7)(2, 5)(3, 8)(4, 6)$, $(1, 4)(6, 7) \in {\rm
Aut}(N_{10})$, we consider the bistellar 2-moves $\kappa_{56}$ and
$\kappa_{57}$ only. Here, $N_{20} = \kappa_{56}(N_{10})$ and
$\kappa_{67}(N_{10}) = N_{16}$.

The removable edges of $N_{11}$ are $14, 24, 56, 57$ and $67$.
Since $(1, 2)(5, 6)(3, 8) \in {\rm Aut}(N_{11})$, we only consider
the bistellar 2-moves $\kappa_{14}$, $\kappa_{56}$ and
$\kappa_{67}$. Here, $N_{22} = \kappa_{14}(N_{11})$,
$\kappa_{56}(N_{11}) = N_{20}$ and $\kappa_{67}(N_{11}) \cong
N_{18}$ (by the map $(2, 4)(5, 7)$).

The removable edges in $N_{12}$ are $12$, $23$, $45$ and $57$.
Here, $\kappa_{12}(N_{12}) \cong N_{22}$ (by the map $(2, 4,
6)$), $\kappa_{23}(N_{12}) = N_{23}$, $\kappa_{45}(N_{12}) \cong
N_{21}$ (by the map $(1, 6, 5, 2, 7, 4)(3, 8)$) and
$\kappa_{57}(N_{12}) \cong N_{18}$ (by the map $(1, 6, 7, 4)$).

The removable edges in $N_{13}$ are $12, 23, 24, 56, 57$ and
$67$. Since $(1, 4)(6, 7) \in {\rm Aut}(N_{13})$, we only consider
$\kappa_{12}$, $\kappa_{23}$, $\kappa_{57}$ and $\kappa_{67}$.
Here, $\kappa_{12}(N_{13}) \cong N_{22}$ (by the map $(2, 7, 5,
4)$),  $\kappa_{23}(N_{13}) = N_{23}$, $\kappa_{57}(N_{13}) \cong
N_{18}$ (by the map $(1, 4)(6, 7)$) and $\kappa_{67}(N_{13}) =
N_{16}$.

The removable edges in $N_{14}$ are $38, 56, 57, 67$. Since $(1,
2, 4)(5, 6, 7)(3, 8) \in {\rm Aut}(N_{14})$, we only consider
$\kappa_{38}$ and $\kappa_{57}$. Here, $N_{24} =
\kappa_{38}(N_{14})$ and $\kappa_{57}(N_{14})= N_{19}$.

The removable edges in $N_{15}$ are $15, 23, 24, 58$. Since
$(1,7)(2,5)(3,8)(4,6) \in {\rm Aut}(N_{15})$, we only consider the
bistellar 2-moves $\kappa_{23}$ and $\kappa_{24}$. Here,
$\kappa_{23}(N_{15}) = N_{23}$ and $\kappa_{24}(N_{15}) \cong
N_{21}$ (by the map $(1, 6, 5, 7, 4)$). This proves the result
for $j=0$.

By the same arguments as in the case for $j=0$, one proves for
the cases for $j = 1, 2$. We summarize these cases in the diagram
given below. Last part follows from the fact that, for $N_{i} \in
{\cal N}_3$, $N_{i}$ has no removable edge. \hfill $\Box$

\smallskip

\setlength{\unitlength}{2.9mm}

\begin{picture}(48.8,20)(1,5)


\thicklines

\put(0,23){\mbox{$N_{10}$}} \put(6,23){\mbox{$N_{7}$}}
\put(12,23){\mbox{$N_{9}$}} \put(18,23){\mbox{$N_{15}$}}
\put(24,23){\mbox{$N_{11}$}} \put(30,23){\mbox{$N_{13}$}}
\put(36,23){\mbox{$N_{12}$}} \put(42,23){\mbox{$N_{8}$}}
\put(48,23){\mbox{$N_{14}$}}

\put(0.5,19.3){\line(0,1){3.2}} \put(1.8,18.9){\line(6,1){22}}

\put(6,19.2){\line(-3,2){4.5}} \put(6.5,19.3){\line(0,1){3.2}}
\put(7,19.2){\line(3,2){4.5}} \put(7.8,19.1){\line(6,1){21.8}}

\put(12.5,19.3){\line(0,1){3.2}} \put(13,19.2){\line(3,2){4.5}}
\put(13.8,19.1){\line(6,1){21.8}}

\put(18,19.2){\line(-3,2){4.5}} \put(19,19.2){\line(0,1){3.4}}
\put(20,19.2){\line(3,1){10}} \put(20.3,18.5){\line(4,1){15.4}}

\put(24.5,19.3){\line(0,1){3.2}} \put(25.6,19.2){\line(5,3){5.2}}
\put(26.5,19){\line(3,1){10}}

\put(30,19.2){\line(-3,2){4.8}} \put(31.5,19.3){\line(0,1){3.2}}
\put(32,19.2){\line(5,3){5.2}} \put(32.7,18.8){\line(5,2){9}}

\put(38,19.2){\line(3,2){4.5}} \put(43,19.3){\line(0,1){3.2}}
\put(43.5,19.2){\line(3,2){4.8}} \put(49,19.3){\line(0,1){3.2}}

\put(0,18){\mbox{$N_{20}$}} \put(6,18){\mbox{$N_{16}$}}
\put(12,18){\mbox{$N_{21}$}} \put(18,18){\mbox{$N_{23}$}}
\put(24,18){\mbox{$N_{22}$}} \put(30,18){\mbox{$N_{18}$}}
\put(36,18){\mbox{$N_{17}$}} \put(42,18){\mbox{$N_{19}$}}
\put(48,18){\mbox{$N_{24}$}}

\put(6,14.2){\line(-3,2){4.5}} \put(6.5,14.3){\line(0,1){3.2}}
\put(7,14.2){\line(3,2){4.7}} \put(7.5,14){\line(5,1){16.5}}

\put(12,14.2){\line(-3,1){10}} \put(12.8,14.2){\line(-3,2){4.5}}
\put(13.5,14.2){\line(5,1){16.5}}

\put(18,14.2){\line(-3,2){4.5}} \put(19,14.2){\line(0,1){3.2}}
\put(20,14.1){\line(3,2){4.7}}

\put(24,14.2){\line(-3,1){10}} \put(26,14.1){\line(3,2){4.7}}

\put(30,14.2){\line(-3,2){4.7}}

\put(36,14.2){\line(-3,1){10}} \put(36.8,14.2){\line(-3,2){4.5}}
\put(37.5,14.2){\line(0,1){3.2}}

\put(42.5,14.2){\line(-3,1){10}} \put(43.5,14.2){\line(0,1){3.2}}

\put(6,13){\mbox{$N_{26}$}} \put(12,13){\mbox{$N_{25}$}}
\put(18,13){\mbox{$N_{30}$}} \put(24,13){\mbox{$N_{29}$}}
\put(30,13){\mbox{$N_{31}$}} \put(36,13){\mbox{$N_{27}$}}
\put(42,13){\mbox{$N_{28}$}}

\put(14.5,9){\line(-2,1){6.5}} \put(15,9.2){\line(-1,2){1.65}}
\put(16.4,9.2){\line(2,1){7.2}} \put(17.5,9){\line(3,1){11.6}}
\put(18,8.7){\line(5,1){18.2}}

\put(21,9){\line(-2,1){7}} \put(23,8.8){\line(5,1){18.8}}

\put(27,9){\line(-2,1){6.5}} \put(27.5,9){\line(-1,2){1.65}}
\put(28,9.3){\line(1,1){3}}

\put(33.2,9.2){\line(-2,1){6.9}} \put(35.3,8.8){\line(2,1){7}}

\put(15,8){\mbox{$N_{32}$}} \put(21,8){\mbox{$N_{33}$}}
\put(27,8){\mbox{$N_{35}$}} \put(33,8){\mbox{$N_{34}$}}

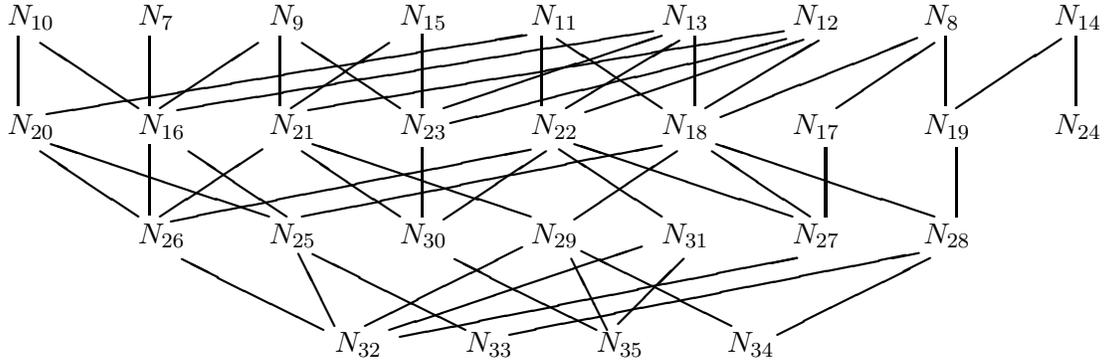
\put(1,5.5){\mbox{{\bf Diagram 2\,:} Hasse diagram of the poset of
all the 3-pseudomanifolds $N_7, \dots, N_{35}$ }}

\end{picture}


\medskip

\noindent {\bf Proof of Corollary \ref{t3}.} Let ${\cal S}_0,
\dots, {\cal S}_6$ be as in the proof of Lemma $\ref{le3.2}$. Let
$M$ be an 8-vertex combinatorial 3-manifold. Then, by Theorem
\ref{t1}, there exist bistellar 1-moves $\kappa_{A_1}, \dots,
\kappa_{A_m}$, for some $m \geq 0$, such that $M_1 :=
\kappa_{A_m}(\cdots (\kappa_{A_1}(M)))$ is a neighbourly 8-vertex
3-pseudomanifold. Since bistellar moves send a combinatorial
3-manifold to a combinatorial 3-manifold, $M_1$ is a
combinatorial 3-manifold. Then, by Theorem \ref{t2}, $M_1 \in
{\cal S}_0$. In other words, $M = \kappa_{e_1}(\cdots
(\kappa_{e_m}(M_1)))$, where $M_1 \in {\cal S}_0$ and
$\kappa_{e_m} : M_1 \mapsto \kappa_{e_m}(M_1)$, $\kappa_{e_i} :
\kappa_{e_{i + 1}}(\cdots (\kappa_{e_m}(M_1)))$ $\mapsto
\kappa_{e_{i}}(\cdots (\kappa_{e_m}(M_1)))$, for $1 \leq i \leq m
- 1$ are bistellar 2-moves. Therefore, by Lemma \ref{S0toS1}, $M
\in {\cal S}_0 \cup \cdots \cup {\cal S}_6$. The result now
follows from Lemma \ref{le3.2}. \hfill $\Box$

\bigskip


\noindent {\bf Proof of Corollary \ref{t4}.} Let ${\cal N}_0,
\dots, {\cal N}_3$ be as in the proof of Lemma $\ref{le3.5}$. Let
$M$ be an 8-vertex normal 3-pseudomanifold. Then, by Theorem
\ref{t1}, there exist bistellar 1-moves $\kappa_{A_1}, \dots,
\kappa_{A_m}$, for some $m\geq 0$, such that $M_1 :=
\kappa_{A_m}(\cdots (\kappa_{A_1}(M)))$ is a neighbourly
3-pseudomanifold. Since bistellar moves send a normal
3-pseudomanifold to a normal 3-pseudomanifold, $M_1$ is normal.
Hence, by Theorem \ref{t2}, $M_1 \in {\cal N}_0$. In other words,
$M = \kappa_{e_1}(\cdots (\kappa_{e_m}(M_1)))$, where $M_1 \in
{\cal N}_0$ and $\kappa_{e_m} : M_1 \mapsto \kappa_{e_m}(M_1)$,
$\kappa_{e_i} : \kappa_{e_{i + 1}}(\cdots(\kappa_{e_m}(M_1)))
\mapsto \kappa_{e_{i}}(\cdots(\kappa_{e_m}(M_1)))$, for $1 \leq i
\leq m-1$ are bistellar 2-moves. Therefore, by Lemma
\ref{N0toN1}, $M \in {\cal N}_0 \cup {\cal N}_1 \cup {\cal N}_2
\cup {\cal N}_3$. The result now follows from Lemma \ref{le3.5}.
\hfill $\Box$

\bigskip

\noindent {\bf Acknowledgement\,:} The authors thank the
anonymous referees for many useful comments which helped to
improve the presentation of this paper. The first author was
partially supported by DST (Grant: SR/S4/MS-272/05) and by
UGC-SAP/DSA-IV.

\end{document}